\documentclass[twoside,11pt,leqno]{article}

\usepackage{amsfonts}
\usepackage{amsmath}

\newtheorem{theorem}{Theorem}[section]
\newtheorem{lemma}[theorem]{Lemma}
\newtheorem{corollary}[theorem]{Corollary}

\long\gdef\boxit#1{\begingroup\vbox{\hrule\hbox{\vrule\kern3pt
      \vbox{\kern3pt#1\kern3pt}\kern3pt\vrule}\hrule}\endgroup}

\def\qed{ \ \vrule width.2cm height.2cm depth0cm\smallskip}

\def\Keywords{\bigskip\par {\sl Keywords\/}:\enspace}

\newenvironment{proof}{\noindent {\bf Proof.\/}}{$\qed$\vskip 0.1in}
\newcommand{\Xcomment}[1]{}

\newcommand{\refeq}[1]{(\ref{eq:#1})}

\newenvironment{myitem}{\refstepcounter{equation}\begin{enumerate}%
\item[(\thesection.\arabic{equation})]$\quad$}{\end{enumerate}}

\makeatletter
\@addtoreset{equation}{section}

\renewcommand{\section}{\@startsection{section}{1}{0pt}%
{-3.5ex plus -1ex minus -.2ex}{2.3ex plus .2ex}%
{\normalfont\Large}}

\renewcommand{\subsection}{\@startsection{subsection}{2}{0pt}%
{-3.25ex plus -1ex minus -.2ex}{1.5ex plus .2ex}%
{\normalfont\large\bf}}

\makeatother

\def\Rset{{\mathbb R}}
\def\Zset{{\mathbb Z}}
\def\Cset{{\mathbb C}}

\def\Bscr{{\cal B}}
\def\Cscr{{\cal C}}

\def\Escr{{\cal E}}
\def\Fscr{{\cal F}}

\def\Hscr{{\cal H}}

\def\Lscr{{\cal L}}
\def\Pscr{{\cal P}}
\def\Rscr{{\cal R}}

\def\Vscr{{\cal V}}

\def\diver{{\rm div}}

\def\tilde{\widetilde}
\def\bar{\overline}
\def\eps{\epsilon}

\def\rest#1{_{\,\vrule height 1.8ex width 0.05em depth 0pt\, #1}}

\def\Vbullet{\Vscr^{\bullet}}
\def\Ebullet{\Escr^{\bullet}}

\begin{document}
\begin{titlepage}

\title{\large\bf Integer Concave Cocirculations and Honeycombs}

\def\thepage {} 

\author{{\large Alexander V. Karzanov}
\\
\ \\
{\sl Institute for System Analysis}\\
{\sl 9, Prospect 60 Let Oktyabrya, 117312 Moscow, Russia}\\
{{\sl E-mail}: {\it sasha@cs.isa.ac.ru}}
}


\maketitle

\begin{abstract}
A {\em convex triangular grid} is represented by a planar digraph
$G$ embedded in the plane so that (a) each bounded face is surrounded
by three edges and forms an equilateral triangle, and (b) the union
$\Rscr$ of bounded faces is a convex polygon.
A real-valued function $h$ on the edges of $G$ is called a {\em
concave cocirculation} if $h(e)=g(v)-g(u)$ for each edge $e=(u,v)$,
where $g$ is a concave function on $\Rscr$ which is affinely linear
within each bounded face of $G$.

Knutson and Tao~\cite{KT} proved an integrality theorem for
so-called {\em honeycombs}, which is equivalent to the assertion that
an integer-valued function on the boundary edges of $G$ is
extendable to an integer concave cocirculation if it is extendable
to a concave cocirculation at all.

In this paper we show a sharper property: for any concave
cocirculation $h$ in $G$, there exists an integer concave
cocirculation $h'$ satisfying $h'(e)=h(e)$ for each boundary edge $e$
with $h(e)$ integer and for each edge $e$ contained in a
bounded face where $h$ takes integer values on all edges.

On the other hand, we explain that for a 3-side grid $G$ of size $n$,
the polytope of concave cocirculations with fixed integer values on
two sides of $G$ can have a vertex $h$ whose entries are integers
on the third side but $h(e)$ has denominator $\Omega(n)$ for some
interior edge $e$.
Also some algorithmic aspects and related results on honeycombs are
discussed.
  \end{abstract}

\medskip
\Keywords Planar graph, Lattice, Discrete convex function, Honeycomb


\end{titlepage}

\baselineskip 15pt

\section{Introduction} \label{sec:intr}

Knutson and Tao~\cite{KT} proved a conjecture concerning
highest weight representations of ${\rm GL}_n(\Cset)$. They
used one combinatorial model, so-called {\em honeycombs}, and an
essential part of the whole proof was to show the existence of an
integer honeycomb under prescribed integer boundary data. The obtained
integrality result for honeycombs admits a re-formulation in terms of
discrete concave functions on triangular grids in the plane.

The purpose of this paper is to show a sharper integrality property
for discrete concave functions.

We start with basic definitions.
Let $\xi_1,\xi_2,\xi_3$ be three affinely independent vectors in the
euclidean plane $\Rset^2$, whose sum is the zero vector.
By a {\em convex (triangular) grid} we mean a finite planar
digraph $G=(V(G),E(G))$ embedded in the plane such that:
(a) each bounded face of $G$ is a triangle surrounded by three edges
and each edge $(u,v)$ satisfies $v-u\in\{\xi_1,\xi_2,\xi_3\}$;
and (b) the region $\Rscr=\Rscr(G)$ of the plane spanned by $G$ is
a convex polygon.
In this paper a convex grid can be considered up to an affine
transformation, and to visualize objects and constructions in what
follows, we will fix the generating vectors $\xi_1,\xi_2,\xi_3$ as
$(1,0)$, $(-1,\sqrt{3})/2$, $(-1,-\sqrt{3})/2$, respectively.
Then each bounded face
is an equilateral triangle ({a \em little triangle} of $G$)
surrounded by a directed circuit with three edges (a {\em
3-circuit}). When $\Rscr$ forms a (big) triangle, we call
$G$ a {\em 3-side grid} (this case is most popular in applications).

A real-valued function $h$ on the edges of $G$ is called a
{\em cocirculation} if the equality $h(e)+h(e')+h(e'')=0$ holds
for each 3-circuit formed by edges $e,e',e''$. This is equivalent to
the existence of a function $g$ on $\Rscr$ which is affinely linear
within each bounded face and satisfies $h(e)=g(v)-g(u)$ for each edge
$e=(u,v)$. Such a $g$ is determined up to adding a constant, and we
refer to $h$ as a {\em concave cocirculation} if $g$ is concave.
(The restriction of such a $g$ to $V(G)$ is usually called a
{\em discrete concave} function.)
It is easy to see that a cocirculation $h$ is concave if
and only if each {\em little rhombus} $\rho$ (the union of two
little triangles sharing a common edge) satisfies the following
{\em rhombus condition}:
 \begin{myitem}
$h(e)\ge h(e')$, where $e,e'$ are non-adjacent (parallel) edges in
$\rho$, and $e$ enters an obtuse vertex of $\rho$.
\label{eq:concave}
  \end{myitem}
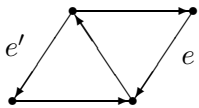
 \begin{center}
  \unitlength=1mm
  \begin{picture}(35,12)
\put(5,0){\circle*{1.0}}
\put(21,0){\circle*{1.0}}
\put(13,12){\circle*{1.0}}
\put(29,12){\circle*{1.0}}
\put(5,0){\vector(1,0){15.5}}
\put(13,12){\vector(1,0){15.5}}
\put(13,12){\vector(-2,-3){7.7}}
\put(29,12){\vector(-2,-3){7.7}}
\put(21,0){\vector(-2,3){7.7}}
\put(4,6){$e'$}
\put(27.5,5){$e$}
  \end{picture}
 \end{center}

Let $B(G)$ denote the set of edges in the boundary of $G$.
Concave cocirculations are closely related (via Fenchel's type
transformations) to honeycombs, and Knutson and Tao's integrality
result on the latter is equivalent to the following.
 \begin{theorem} {\rm \cite{KT}} \label{tm:KT}
For a convex grid $G$ and a function $h_0:B(G)\to\Zset$, there exists
an integer concave cocirculation in $G$ coinciding with $h_0$ on
$B(G)$ if $h_0$ is extendable to a concave cocirculation in $G$ at all.
  \end{theorem}

For a direct proof of this theorem for 3-side grids, without appealing
to honeycombs, see Buch~\cite{Bu}.
(Note also that for a 3-side grid $G$, a combinatorial
characterization for the set of functions $h_0$ on $B(G)$ extendable
to concave cocirculations is given in~\cite{KTW}: it is a polyhedral
cone in $\Rset^{B(G)}$ whose nontrivial facets are described by
Horn's type inequalities with respect to so-called {\em puzzles}; for
an alternative proof and an extension to arbitrary convex grids,
see~\cite{Kar}.)

In this paper we extend Theorem~\ref{tm:KT} as follows.
\begin{theorem} \label{tm:main}
Let $h$ be a concave cocirculation in a convex grid $G$. There exists
an integer concave cocirculation $h'$ in $G$ such that $h'(e)=h(e)$
for all edges $e\in O_h\cup I_h$.
Here $O_h$ is the set of boundary edges $e$ where $h(e)$ is an integer,
and $I_h$ is the set of edges contained in little triangles $\Delta$
such that $h$ takes integer values on the three edges of $\Delta$.
  \end{theorem}

\noindent {\bf Remark 1.} One could attempt to further strengthen
Theorem~\ref{tm:KT} by asking: can one improve any concave
cocirculation $h$ to an integer concave cocirculation preserving the
values on {\em all} edges where $h$ is integral? In general, the
answer is negative; a counterexample will be given in the end of this
paper.

\medskip
Our method of proof of Theorem~\ref{tm:main} is constructive and
based on iteratively transforming the current concave cocirculation
until the desired integer concave cocirculation is found.
As a consequence, we obtain a polynomial-time combinatorial algorithm
to improve $h$ to $h'$ as required. (The idea of proof of
Theorem~\ref{tm:KT} in~\cite{Bu} is to show the existence of a concave
cocirculation coinciding with $h_0$ on $B(G)$ whose values are
expressed as integer combinations of values of $h_0$; \cite{KT}
establishes an analogous property for honeycombs. Our approach is
different.) We prefer to describe an iteration by considering a
corresponding task on the related honeycomb model and then translating
the output to the language of cocirculations in $G$, as this makes our
description technically simpler and more enlightening.

\smallskip
The above theorems admit a re-formulation in polyhedral terms.
Given a subset $F\subseteq E(G)$ and a function $h_0:F\to \Rset$,
let $\Cscr(G,h_0)$ denote the set of concave cocirculations in $G$
such that $h(e)=h_0(e)$ for all $e\in F$.
Since concave cocirculations are described by linear constraints,
$\Cscr(G,h_0)$ forms a (possibly empty) polyhedron in $\Rset^{E(G)}$.
Then Theorem~\ref{tm:KT} says that such a polyhedron (if nonempty) has
an integer point $h$ in the case $h_0$ is an integer-valued function
on $B(G)$, whereas Theorem~\ref{tm:main} is equivalent to saying
that a similar property takes place if $h_0:F\to\Zset$ and
$F=B'\cup\cup(E(\Delta):\Delta\in T)$, where $B'\subseteq B(G)$ and
$T$ is a set of little triangles. (Note that when $F=B(G)$ and when
$\Rscr(G)$ is not a hexagon, one can conclude from the concavity that
$\Cscr(G,h_0)$ is bounded, i.e., it is a polytope.)

\smallskip
On the ``negative'' side, it turned out that $\Cscr(G,h_0)$ with
$h_0:B(G)\to \Zset$ need not be an integer polytope; an example
with a half-integer but not integer vertex is given in~\cite{KT}
(and in~\cite{Bu}). One can show that the class of such polyhedra
has ``unbounded fractionality''. Moreover, denominators of vertex
entries can be arbitrarily increasing as the size of $G$ grows even if
functions $h_0$ with smaller domains are considered. Hereinafter by
the {\em size} of $G$ we mean its maximum side length (=number of
edges). We show the following.
\begin{theorem} \label{tm:frac}
For any positive integer $k$, there exists a 3-side grid $G$ of size
$O(k)$ and a function $h_0:F\to\Zset$, where $F$ is the set of
edges of two sides of $G$, such that the polyhedron $\Cscr(G,h_0)$ has
a vertex $h$ satisfying: (a) $h(e)$ has denominator $k$ for some edge
$e\in E(G)$, and (b) $h$ takes integer values on all boundary edges.
  \end{theorem}

(One can see that in this case $\Cscr(G,h_0)$ is also a polytope.
Note also that if $h'_0$ is the restriction of $h$ to a set
$F'\supset F$, then $h$ is, obviously, a vertex of the polytope
$\Cscr(G,h'_0)$ as well.)

This paper is organized as follows.
In Section~\ref{sec:honey} we explain the notion of
honeycomb and a relationship between honeycombs and concave
cocirculations. Sections~\ref{sec:l_path} and~\ref{sec:deform}
consider special paths (open and closed) in a honeycomb and describe
a certain transformation of the honeycomb in a
neighbourhood of such a path. Section~\ref{sec:proof} explains how
to ``improve'' the honeycomb by use of such transformations and
eventually proves Theorem~\ref{tm:main}. A construction of $G,h_0$
proving Theorem~\ref{tm:frac} is given in Section~\ref{sec:frac};
it relies on an approach involving honeycombs as well. We also explain
there (in Remark 3) that the set $F$ in this theorem can be reduced
further. The concluding Section~\ref{sec:concl} discusses algorithmic
aspects, suggests a slight strengthening of Theorem~\ref{tm:main},
gives a counterexample mentioned in Remark 1, and raises an open
question.

\section{Honeycombs}
\label{sec:honey}

For technical needs of this paper, our definition of honeycombs
will be somewhat different from, though equivalent to, that given
in~\cite{KT}. It is based on a notion of pre-honeycombs, and before
introducing the latter, we clarify some terminology and notation user
later on. Let $\xi_1,\xi_2,\xi_3$ be the generating vectors as above
(note that they follow anticlockwise around the origin).

The term {\em line} is applied to {\em (fully) infinite},
{\em semiinfinite}, and {\em finite} lines, i.e., to sets of the form
$a+\Rset b$, $a+\Rset_+ b$, and $\{a+\lambda b: 0\le\lambda\le 1\}$,
respectively, where $a\in\Rset^2$ and
$b\in\Rset^2\setminus\{{\bf 0}\}$.
For a vector (point) $v\in\Rset^2$ and $i=1,2,3$, we denote by
$\Xi_i(v)$ the infinite line containing $v$ and perpendicular to
$\xi_i$, i.e., the set of points $u$ with $(u-v)\cdot \xi_i=0$.
(Hereinafter $x\cdot y$ denotes the inner product of vectors $x,y$.)
The line $\Xi_i(v)$ is the union of two semiinfinite lines
$\Xi_i^+(v)$ and $\Xi_i^-(v)$ with the end $v$, where
the rays $\Xi_i^+(v)-v$, $\Rset_+\xi_i$ and $\Xi_i^-(v)-v$
follow in the {\em anticlockwise} order around the origin.
Any line perpendicular to $\xi_i$ is called a $\Xi_i$-{\em line}.

By a $\Xi$-{\em system} we mean a finite set $\Lscr$ of $\Xi_i$-lines
($i\in\{1,2,3\}$) along with an {\em integer} weighting $w$ on them.
For a point $v\in\Rset^2$, a ``sign'' $s\in\{+,-\}$, and $i=1,2,3$,
define $w^s_i(v)$ to be the sum of weights $w(L)$ of the
lines $L\in\Lscr$ whose intersection with $\Xi^s_i(v)$ contains $v$ and
is a line (not a point). We call a $\Xi$-system $(\Lscr,w)$ a
{\em pre-honeycomb} if for any point $v$, the numbers
$w^s_i(v)$ are {\em nonnegative} and satisfy the condition
  \begin{equation}  \label{eq:div}
  w_1^+(v)-w_1^-(v)=w_2^+(v)-w_2^-(v)=w_3^+(v)-w_3^-(v)=:\diver_w(v);
  \end{equation}
$\diver_w(v)$ is called the {\em divergency} at $v$.

\begin{figure}[tb]
 \begin{center}
  \unitlength=1mm
  \begin{picture}(150,36)
\put(8,18){\circle*{1.0}}
\put(8,18){\vector(1,0){16}}
\put(8,18){\vector(-2,-3){8}}
\put(8,18){\vector(-2,3){8}}
\put(17,14){$\xi_1$}
\put(4,26){$\xi_2$}
\put(4,8){$\xi_3$}
\put(68,18){\circle*{1.0}}
\put(68,18){\line(0,-1){18}}
\put(68,18){\line(3,2){15}}
\put(68,18){\line(-3,2){15}}
\put(67,20){$v$}
\put(69.5,4){$\Xi^+_1(v)$}
\put(51,20){$\Xi^+_3(v)$}
\put(78,20){$\Xi^+_2(v)$}
\put(135,18){\circle*{1.0}}
\put(135,18){\line(0,-1){18}}
\put(135,18){\line(3,-2){15}}
\put(135,18){\line(3,2){15}}
\put(135,18){\line(-3,2){24}}
\put(135,18){\line(-3,-2){24}}
\put(123,0){\line(0,1){36}}
\put(123,10){\circle*{1.0}}
\put(123,26){\circle*{1.0}}
\put(116,7){1}
\put(116,27){3}
\put(120.5,17){1}
\put(128,11){1}
\put(128,23){3}
\put(143,9){1}
\put(143,25){3}
\put(124,2){1}
\put(124,33){1}
\put(136,4){2}
\put(138,17){$v$}
\put(124.5,8){$u$}
\put(124.5,27){$z$}
  \end{picture}
 \end{center}
\caption{Generating vectors, lines $\Xi^+_i(v)$, and a honeycomb
instance} \label{fig:hon}
  \end{figure}
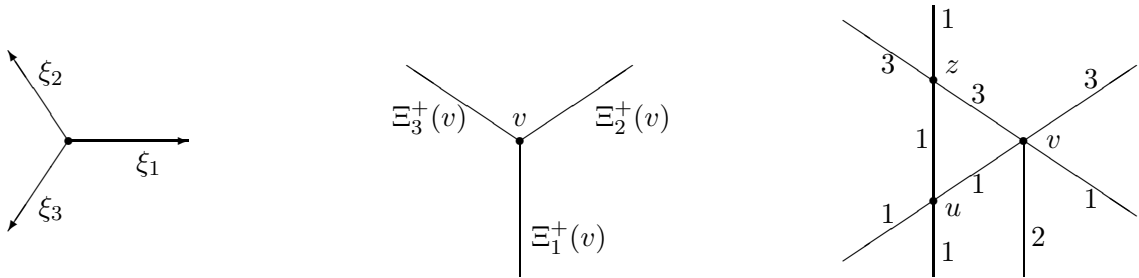

Now a {\em honeycomb} is a certain non-standard edge-weighted planar
graph $\Hscr=(\Vscr,\Escr,w)$ with vertex set $\Vscr\ne\emptyset$
and edge set $\Escr$ in which non-finite edges are allowed.
More precisely: (i) each vertex is incident with at least 3 edges;
(ii) each edge is a line with no interior point contained in another
edge; (iii) $w(e)$ is a {\em positive} integer for each $e\in\Escr$;
and (iv) $(\Escr,w)$ is a pre-honeycomb.
Then each vertex $v$ has degree at most 6, and for $i=1,2,3$ and
$s\in\{+,-\}$, we denote by $e_i^s(v)$ the edge incident to $v$ and
contained in $\Xi^s_i(v)$ when such an edge exists, and say that this
edge {\em has sign} $s$ at $v$.
In~\cite{KT} the condition on a vertex $v$ of a honeycomb similar
to~\refeq{div} is called the {\em zero-tension condition}, motivated by
the observation that if each edge incident to $v$ pulls on $v$ with a
tension equal to its weight, then the total force applied to $v$ is
zero.
Figure~\ref{fig:hon} illustrates a honeycomb with three vertices
$u,v,z$ and ten edges of which seven are semiinfinite.

We will take advantage of the fact that any pre-honeycomb
$(\Lscr,w)$ determines, in a natural way, a unique honeycomb
$\Hscr=(\Vscr,\Escr,w')$ with $w^s_i(v)=(w')^s_i(v)$ for all $v,s,i$.
Here $\Vscr$ is the set of points $v$ for which at least three
numbers among $w^s_i(v)$'s are nonzero. The set $\Escr$ consists of
all maximal $\Xi_i$-lines $e$ for $i=1,2,3$ such that any interior
point $v$ on $e$ satisfies $w^+_i(v)>0$ and does not belong to
$\Vscr$; the weight $w'(e)$ is defined to be just this number
$w^+_i(v)$ (=$w^-_i(v)$), which does not depend on $v$.

Since $\Vscr\ne\emptyset$, one can conclude from~\refeq{div} that
a honeycomb has no fully infinite edge but the set of
semiinfinite edges in it is nonempty. This set, called the
{\em boundary} of $\Hscr$ and denoted by $\Bscr(\Hscr)$, is naturally
partitioned into subsets $\Bscr_i^s$ consisting of the
semiinfinite edges of the form $\Xi^s_i(\cdot)$.
Then~\refeq{div} implies that $w(\Bscr^+_i)-w(\Bscr^-_i)$ is the
same for $i=1,2,3$.
(For a subset $E'\subseteq E$ and a function
$c:E\to\Rset$, $c(E')$ stands for $\sum(c(e):e\in E')$.)

Let us introduce the {\em dual coordinates} $d_1,d_2,d_3$ of a point
$x\in\Rset^2$ by setting
  $$
d_i(x):=-x\cdot\xi_i, \qquad i=1,2,3.
  $$
Since $\xi_1+\xi_2+\xi_3=0$, one has
  \begin{equation} \label{eq:zero_sum}
      d_1(x)+d_2(x)+d_3(x)=0\qquad \mbox{for each $x\in\Rset^2$}.
  \end{equation}
When one traverses an edge $e$ of a honeycomb, one dual coordinate
remains constant while the other two trade off; this {\em constant
dual coordinate} is denoted by $d^c(e)$.

Next we explain that the honeycombs one-to-one correspond
to the concave cocirculations via a sort of planar duality.
Consider a honeycomb $\Hscr=(\Vscr,\Escr,w)$. Let $v\in\Vscr$.
Since the numbers $w^s_i(v)$ are nonnegative,
condition~\refeq{div} is equivalent to the existence of a (unique)
grid $G_v$ whose boundary is formed, in the anticlockwise
order, by $w_1^+(v)$ edges parallel to $\xi_1$, followed by
$w_3^-(v)$ edges parallel to $\xi_3$, followed by $w_2^+(v)$ edges
parallel to $\xi_2$, and so on, as illustrated in the picture.
 \begin{center}
  \unitlength=1mm
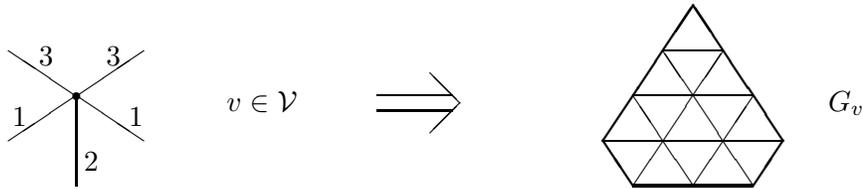
  \begin{picture}(120,24)
\put(10,12){\circle*{1.0}}
\put(10,12){\line(0,-1){12}}
\put(10,12){\line(3,-2){9}}
\put(10,12){\line(3,2){9}}
\put(10,12){\line(-3,2){9}}
\put(10,12){\line(-3,-2){9}}
\put(11,2){2}
\put(1.5,8){1}
\put(17,8){1}
\put(5,16){3}
\put(14,16){3}
\put(30,10){$v\in\Vscr$}
\put(50,12){\line(1,0){10}}
\put(50,10){\line(1,0){10}}
\put(61,11){\line(-1,1){4}}
\put(61,11){\line(-1,-1){4}}
\put(80,6){\line(1,0){24}}
\put(84,12){\line(1,0){16}}
\put(88,18){\line(1,0){8}}
\put(84,0){\line(2,3){12}}
\put(92,0){\line(2,3){8}}
\put(92,0){\line(-2,3){8}}
\put(100,0){\line(-2,3){12}}
\put(110,10){$G_v$}
\thicklines
\put(84,0){\line(1,0){16}}
\put(80,6){\line(2,3){12}}
\put(100,0){\line(2,3){4}}
\put(84,0){\line(-2,3){4}}
\put(104,6){\line(-2,3){12}}
  \end{picture}
 \end{center}
For an edge $e$ incident to $v$, label $e^\ast$ the corresponding
side of $G_v$. For each finite edge
$e=uv$ of $\Hscr$, glue together the grids $G_u$ and $G_v$ by
identifying the sides labelled $e^\ast$ in both. One can see
that the resulting graph $G$ is a convex grid. Also each nonempty
set $\Bscr_i^s$ one-to-one corresponds to a side of $G$, denoted by
$B_i^s$; it is formed by $w(\Bscr_i^s)$ edges parallel to $\xi_i$,
and the outward normal at $B_i^s$ points at the direction $(i,s)$.
The picture below illustrates the grid generated by the honeycomb in
Fig.~\ref{fig:hon}.
 \begin{center}
  \unitlength=1mm
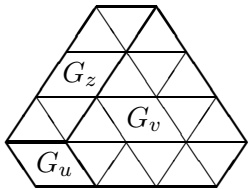
  \begin{picture}(40,24)
\put(8,6){\line(1,0){24}}
\put(4,12){\line(1,0){24}}
\put(8,18){\line(1,0){16}}
\put(12,0){\line(2,3){4}}
\put(20,12){\line(2,3){4}}
\put(20,0){\line(2,3){8}}
\put(8,6){\line(-2,3){4}}
\put(20,0){\line(-2,3){8}}
\put(28,0){\line(-2,3){16}}
\put(4,2){$G_u$}
\put(16,8){$G_v$}
\put(7.5,14){$G_z$}
\thicklines
\put(4,0){\line(1,0){24}}
\put(0,6){\line(1,0){8}}
\put(12,24){\line(1,0){8}}
\put(0,6){\line(2,3){12}}
\put(8,6){\line(2,3){12}}
\put(28,0){\line(2,3){4}}
\put(4,0){\line(-2,3){4}}
\put(12,0){\line(-2,3){4}}
\put(32,6){\line(-2,3){12}}
  \end{picture}
 \end{center}
The dual coordinates of vertices of $\Hscr$ generate the function
$h$ on $E(G)$ defined by:
  \begin{equation}  \label{eq:h-d}
  h(e):=d_i(v), \quad
  \mbox{where $v\in\Vscr$, $i=1,2,3$, and $e$ is an edge in $G_v$
parallel to $\xi_i$}.
  \end{equation}

Then~\refeq{zero_sum} implies $h(e)+h(e')+h(e'')=0$ for any little
triangle with edges $e,e',e''$ in $G$, i.e., $h$ is a cocirculation.
To see that $h$ is concave, it suffices to check~\refeq{concave} for a
little rhombus $\rho$ formed by little triangles lying in different
graphs $G_u$ and $G_v$.
Then $\tilde e=uv$ is an edge of $\Hscr$; let $\tilde e$ be
perpendicular to $\xi_i$ and assume $\tilde e=e_i^-(u)=e_i^+(v)$.
Observe that $d_{i+1}(u)>d_{i+1}(v)$ (taking indices modulo 3) and
that the side-edge $e$ of $\rho$ lying in $G_u$ and parallel to
$\xi_{i+1}$ enters an obtuse vertex of $\rho$. Therefore,
$h(e)=d_{i+1}(u)>d_{i+1}(v)=h(e')$, where $e'$ is the side-edge of
$\rho$ parallel to $e$ (lying in $G(v)$), as required.

Conversely, let $h$ be a concave cocirculation in a convex grid $G$.
Subdivide $G$ into maximal subgraphs $G_1,\ldots,G_k$, each being the
union of little triangles where, for each $i=1,2,3$, all edges parallel
to $\xi_i$ have the same value of $h$. The concavity of $h$ implies
that each $G_j$ is again a convex grid; it spans a maximal region
where the corresponding function $g$ on $\Rscr$ is affinely linear,
called a {\em flatspace} of $h$.
For $j=1,\ldots,k$, take the point $v_j$ in the
plane defined by the dual coordinates $d_i(v_j)=h(e_i)$, $i=1,2,3$,
where $e_i$ is an edge of $G_j$ parallel to $\xi_i$.
(The property $h(e_1)+h(e_2)+h(e_3)=0$ implies~\refeq{zero_sum},
so $v_j$ exists; also the points $v_j$ are different).
For each pair of graphs $G_j,G_{j'}$ having a common side
$S$, connect $v_j$ and $v_{j'}$ by line (finite edge) $\ell$; observe
that $\ell$ is perpendicular to $S$.
And if a graph $G_j$ has a side $S$ contained in the boundary of
$G$, assign the semiinfinite edge $\ell=\Xi_i^s(v_j)$ whose direction
$(i,s)$ corresponds to the outward normal at $S$. In all
cases the weight of $\ell$ is assigned to be the number of edges in
$S$. One can check (by reversing the argument above) that the obtained
sets of points and weighted lines constitute a honeycomb $\Hscr$, and
that the above construction for $\Hscr$ returns $G,h$.

\section{Legal Paths and Cycles}
\label{sec:l_path}


In this section we consider certain paths (possibly cycles) in
a honeycomb $\Hscr=(\Vscr,\Escr,w)$.
A transformation of $\Hscr$ with respect to such a path,
described in the next section, ``improves'' the honeycomb, in a
certain sense, and we will show that a number of such improvements
results in a honeycomb determining an integer concave cocirculation as
required in Theorem~\ref{tm:main}. First of all we need some
definitions and notation.

Let $\Vbullet$ denote the set of vertices $v\in\Vscr$ having at least
one nonintegral dual coordinate $d_i(v)$, and $\Ebullet$ the set
of edges $e\in\Escr$ whose constant coordinate $d^c(e)$ is
nonintegral. For brevity we call such vertices and edges
{\em nonintegral}.

For $s\in\{+,-\}$, $-s$ denotes the sign opposite to $s$.
Let $v\in\Vscr$.
An edge $e_i^s(v)$ is called {\em dominating} at $v$ if
$w_i^s(v)>w_i^{-s}(v)$. By~\refeq{div}, $v$ has either none or three
dominating edges, each pair forming an angle of $120^\circ$.
A pair $\{e_i^s(v),e_j^{s'}(v)\}$ of distinct {\em nonintegral} edges
is called {\em legal} if either they are {\em opposite} to each
other at $v$, i.e., $j=i$ and $s'=-s$, or both edges are
dominating at $v$ (then $j\ne i$ and $s'=s$).
By a {\em path} in $\Hscr$ we mean a finite alternating sequence
$P=(v_0,q_1,v_1,\ldots,q_k,v_k)$, $k\ge 1$, of vertices and edges where:
for $i=2,\ldots,k-1$, $q_i$ is
a finite edge and $v_{i-1},v_i$ are its ends; $q_1$ is either a finite
edge with the ends $v_0,v_1$, or a semiinfinite edge with the end
$v_1$; similarly, if $k>1$ then $q_k$ is either a finite edge with
the ends $v_{k-1},v_k$, or a semiinfinite edge with the end $v_{k-1}$.
When $q_1$ is semiinfinite, $v_0$ is thought of as a {\em dummy}
(``infinite'') vertex, and we write $v_0=\{\emptyset\}$; similarly,
$v_k=\{\emptyset\}$ when $q_k$ is semiinfinite and $k>1$.
Self-intersecting paths are admitted. We call $P$

(i) an {\em open} path if $k>1$ and both edges $q_1,q_k$ are
semiinfinite;

(ii) a {\em legal path} if each pair of consecutive edges
$q_i,q_{i+1}$ is legal;

(iii) a {\em legal cycle} if it is a legal path,
$v_0=v_k\ne\{\emptyset\}$ and $\{q_k,q_1\}$ is a legal pair.

A legal cycle $P$ is usually considered up to shifting cyclically and
the indices are taken modulo $k$. We say that a legal path $P$
{\em turns} at $v\in\Vscr$ if, for some $i$ with $v_i=v$,
the (existing) edges $q_i,q_{i+1}$ are not opposite at $v_i$ (then
$q_i,q_{i+1}$ are different dominating edges at $v_i$). We also call
such a triple $(q_i,v_i,q_{i+1})$ a {\em bend} of $P$ at $v$.

Assume $\Vbullet$ is nonempty. (When $\Vbullet=\emptyset$, the concave
cocirculation in $G$ determined by $\Hscr$ is already integral.) A
trivial but important observation from~\refeq{zero_sum} is that if a
vertex $v$ has a nonintegral dual coordinate, then it has at least two
nonintegral dual coordinates. This implies $\Ebullet\ne\emptyset$.
Moreover, if $e$ is a nonintegral edge dominating at $v$, then $e$
forms a legal pair with another nonintegral edge dominating at $v$.

Our method of proof of Theorem~\ref{tm:main} will rely on the
existence of a legal path with some additional properties, as follows.

  \begin{lemma}  \label{lm:lpath}
There exists an open legal path or a legal cycle $P=(v_0,q_1,v_1,
\ldots,q_k,v_k)$ such that:

{\rm (i)} each edge $e$ of $\Hscr$ occurs in $P$ at most twice, and
if it occurs exactly twice, then $P$ traverses $e$ in both directions,
i.e., $e=q_i=q_j$ and $i<j$ imply $v_i=v_{j-1}$;

{\rm (ii)} if an edge $e$ occurs in $P$ twice, then $w(e)>1$;

{\rm (iii)} for each vertex $v$ of $\Hscr$, the number of times $P$ turns
at $v$ does not exceed $\min\{2,|\diver_w(v)|\}$.
  \end{lemma}
  \Xcomment{ For the abstract:
(The proof (omitted here) relies on a procedure to construct such a
$P$.)
   }
  \begin{proof}
We grow a legal path $P$, step by step, by the following process.
Initially, choose a nonintegral semiinfinite edge $e$ if it exists, and
set $P=(v_0,q_1,v_1)$, where $v_0:=\{\emptyset\}$, $q_1:=e$, and
$v_1$ is the end of $e$.
Otherwise we start with $P=(v_0,q_1,v_1)$, where $q_1$ is an arbitrary
nonintegral finite edge and $v_0,v_1$ are its ends. Let
$P=(v_0,q_1,v_1\ldots,q_i,v_i)$ be a current legal path with
$v:=v_i\in\Vscr$ satisfying (i),(ii),(iii). At an iteration, we wish
either to increase $P$ by adding some $q_{i+1},v_{i+1}$ (maintaining
(i),(ii),(iii)) or to form the desired cycle.

By the above observation, $e:=q_i$ forms a legal pair with at least
one edge $e'$ incident to $v$. We select such an $e'$ by rules
specified later and act as follows. Suppose $e'$ occurs in $P$ and is
traversed from $v$, i.e., $v=v_{j-1}$ and $e'=e_j$ for some $j<i$.
Then the part of $P$ from $v_{j-1}$ to $v_i$ forms a legal cycle; we
finish the process and output this cycle. Clearly it satisfies (i)
and (ii) (but the number of bends at $v$ may increase). Now suppose
$e'$ is not traversed from $v$. Then we grow $P$ by adding $e'$ as the
new last edge $q_{i+1}$ and adding $v_{i+1}$ to be the end of $e'$
different from $v$ if $e'$ is finite, and to be $\{\emptyset\}$ if
$e'$ is semiinfinite. In the latter case, the new $P$ is an open legal
path (taking into account the choice of the first edge $q_1$); we
finish and output this $P$. And in the former case, we continue the
process with the new current $P$. Clearly property (i) is maintained.

We have to show that $e'$ can be chosen so as to maintain the
remaining properties (concerning $e'$ and $v$). Consider two cases.

\medskip
{\em Case 1.} $e$ is not dominating at $v$. Then $e'$ is opposite to
$e$ at $v$ (as the choice is unique), and (iii) remains valid as no
new bend at $v$ arises. If $e'$ is dominating at $v$, then $w(e')>
w(e)\ge 1$, implying (ii). And if $e'$ is not dominating at $v$, then
$w(e')=w(e)$ and, obviously, the new $P$ traverses $e'$ as many times
as it traverses $e$, implying (ii) as well.

\medskip
{\em Case 2.} $e$ is dominating at $v$. Let the old $P$ turn $b$
times at $v$. First suppose $P$ has a bend $\beta=(q_j,v_j,q_{j+1})$ at
$v$ not using the edge $e$. Since the edges occurring in any bend are
nonintegral and dominating at the corresponding vertex,
$\{e,q_{j+1}\}$ is a legal pair. We choose $e'$ to be $q_{j+1}$. This
leads to forming a cycle with $b$ bends at $v$ as before (as the bend
$\beta$ is destroyed while the only bend $(e,v,e')$ is added),
implying (iii).

So assume $\beta$ as above does not exist. Then the old $P$ can have
at most one bend at $v$, namely, one of the form $\beta'=(q,v,e)$,
whence $b\le 1$. If $b<|\diver_w(v)|$, then taking as $e'$ a
nonintegral edge dominating at $v$ and different from $e$ maintains
both (ii) and (iii) (to see (ii), observe that the number of times $P$
traverses $e'$ is less than $w(e')$). Now let $b=|\diver_w(v)|$. Then
$b=1$ (as $\diver_w(v)\ne 0$) and $P$ has the bend $\beta'$ as above.
Therefore, $P$ traverses $e$ twice (in $\beta'$ and as $q_i$), and we
conclude from this fact together with $|\diver_w(v)|=1$ that $e$ has
the opposite edge $\bar e$ at $v$. Moreover, $\bar e$ cannot occur in
$P$. For otherwise $P$ would traverse $e$ more than twice, taking into
account that $\bar e$ forms a legal pair only with $e$ (as $\bar e$ is
non-dominating at $v$). Thus, the choice of $e'$ to be $\bar e$
maintains (ii) and (iii), completing the proof of the lemma.
  \end{proof}

\section{$\eps$-Deformation}
\label{sec:deform}

Let $P=(v_0,q_1,v_1,\ldots,q_k,v_k)$ be as in Lemma~\ref{lm:lpath}.
Our transformation of the honeycomb $\Hscr=(\Vscr,\Escr,w)$ in
question is, roughly speaking, a result of ``moving a unit weight copy
of $P$ in a normal direction'' (considering $P$ as a curve in the
plane); this is analogous to an operation on more elementary paths or
cycles of honeycombs in~\cite{KT,KTW}. It is
technically easier for us to describe such a transformation by
handling a pre-honeycomb behind $\Hscr$ in which the line set include
the maximal straight subpaths of $P$.

When $(q_i,v_i,q_{i+1})$ is a bend, we say that $v_i$ is a {\em bend
vertex} of $P$. We assume that $v_0$ is a bend vertex if $P$ is
a cycle. For a bend vertex $v_i$, we will distinguish between the
cases when $P$ {\em turns right} and {\em turns left} at $v_i$,
defined in a natural way regarding the orientation of $P$.
Let $v_{t(0)},v_{t(1)},\ldots,v_{t(r)}$ ($0=t(0)<t(1)<\ldots<t(r)=k$)
be the sequence of bend or dummy vertices of $P$. Then for
$i=1,\ldots,r$, the union of edges $q_{t(i-1)+1},\ldots,q_{t(i)}$ is
a (finite, semiinfinite or even fully infinite) line, denoted by $L_i$.
For brevity $v_{t(i)}$ is denoted by $u_i$.

Our transformation of $\Hscr$ depends on a real parameter $\eps>0$
measuring the distance of moving $P$ and on a direction of moving;
let for definiteness we wish to move $P$ ``to the right'' (moving
``to the left'' is symmetric). We assume that $\eps$ is small enough;
an upper bound on $\eps$ will be discussed later.
By the transformation, a unit weight copy of each line $L_i$
(considered as oriented from $u_{i-1}$ to $u_i$) is split off the
honeycomb and moves (possibly extending or shrinking) at distance
$\eps$ to the right, turning into a parallel line $L'_i$ connecting
$u'_{i-1}$ and $u'_i$. Let us describe this construction more formally.
First, for a bend vertex $u_i$, let
the constant coordinates of the lines $L_i$ and $L_{i+1}$ be $p$-th
and $p'$-th dual coordinates, respectively. Then the point $u'_i$ is
defined by
  \begin{myitem}
$d_p(u'_i):=d_p(u_i)-\eps$ and $d_{p'}(u'_i):=d_{p'}(u_i)+\eps$ if both
$L_i,L_{i+1}$ have sign $+$ at $v$, and $d_p(u'_i):=d_p(u_i)+\eps$
and $d_{p'}(u'_i):=d_{p'}(u_i)-\eps$ if $L_i,L_{i+1}$ have sign $-$ at
$v$,
   \label{eq:shiftP}
  \end{myitem}
where, similar to the edges, a $\Xi_i$-line is said to have sign $s$
at its end $v$ if it is contained in $\Xi^s_i(v)$. Possible cases are
illustrated in the picture.
 \begin{center}
  \unitlength=1mm
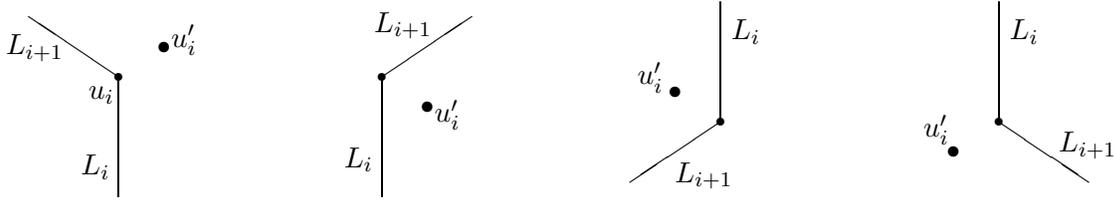
  \begin{picture}(152,26)
  \put(0,0){\begin{picture}(32,26)
\put(15,16){\circle*{1.0}}
\put(15,16){\line(-3,2){12}}
\put(15,16){\line(0,-1){16}}
\put(21,20){\circle*{1.5}}
\put(11,13){$u_i$}
\put(10,3){$L_i$}
\put(0,19){$L_{i+1}$}
\put(22,20){$u'_i$}
\end{picture}}
  \put(40,0){\begin{picture}(32,26)
\put(10,16){\circle*{1.0}}
\put(10,16){\line(3,2){12}}
\put(10,16){\line(0,-1){16}}
\put(16,12){\circle*{1.5}}
\put(5,4){$L_i$}
\put(9,22){$L_{i+1}$}
\put(17,10){$u'_i$}
\end{picture}}
  \put(80,0){\begin{picture}(32,26)
\put(15,10){\circle*{1.0}}
\put(15,10){\line(-3,-2){12}}
\put(15,10){\line(0,1){16}}
\put(9,14){\circle*{1.5}}
\put(16.5,21){$L_i$}
\put(9,2){$L_{i+1}$}
\put(4,15){$u'_i$}
\end{picture}}
  \put(120,0){\begin{picture}(32,26)
\put(12,10){\circle*{1.0}}
\put(12,10){\line(3,-2){12}}
\put(12,10){\line(0,1){16}}
\put(6,6){\circle*{1.5}}
\put(13.5,21){$L_i$}
\put(20,6){$L_{i+1}$}
\put(2,8){$u'_i$}
\end{picture}}
  \end{picture}
 \end{center}

Second, for $i=1,\ldots,r$, define $L'_i$ to be the line connecting
$u'_{i-1}$ and $u'_i$ (when $P$ is an open path, $u'_0$ and $u'_r$ are
dummy points and the non-finite lines $L'_1,L'_r$ are defined in a
natural way). Denoting by $\ell(L)$ the euclidean length
(scalled by $2/\sqrt{3}$) of a line $L$, one can see that if a line
$L_i$ is finite and $\eps\le\ell(L_i)$, then
  \begin{myitem}
$\ell(L'_i)=\ell(L_i)-\eps$ if $P$ turns right at both
$u_{i-1}$ and $u_i$, and $\ell(L'_i)\ge\ell(L_i)$ otherwise.
   \label{eq:length}
  \end{myitem}
This motivates a reasonable upper bound on $\eps$, to be
the minimum length $\bar\eps_0$ of an $L_i$ such that
$P$ turns right at both $u_{i-1},u_i$ ($\bar\eps_0=\infty$ when no
such $L_i$ exists).

Third, consider the $\Xi$-system $\Pscr=
(\{L_1,\ldots,L_r\}\cup\Escr,\tilde w)$, where $\tilde w(L_i)=1$
for $i=1,\ldots,r$, and $\tilde w(e)$ is equal to $w(e)$ minus the
number of occurrencies of $e\in\Escr$ in $L_1,\ldots,L_r$.
This $\Pscr$ is a pre-honeycomb representing $\Hscr$, i.e., satisfying
$w^s_i(v)={\tilde w}^s_i(v)$ for all $v,i,s$.
We replace in $\Pscr$ the lines $L_1,\ldots,L_r$ by the lines
$L'_1,\ldots,L'_r$ with unit weight each.
(When $\eps=\bar\eps_0<\infty$, the length of at least one line $L'_i$
reduces to zero, by~\refeq{length}, and this line vanishes in
$\Pscr'$.) Also for each bend vertex $u_i$, we add line
$R_i$ connecting $u_i$ and $u'_i$. We assign to $R_i$ weight
1 if $P$ turns right at $u_i$, and $-1$ otherwise.
Let $\Pscr'=(\Lscr',w')$ be the resulting $\Xi$-system.

The transformation in a neighbourhood of a bend vertex $v=u_i$
is illustrated in the picture; here the numbers on edges indicate
their original weights or the changes due to the transformation, and
(for simplicity) $P$ passes $v$ only once.
 \begin{center}
  \unitlength=1mm
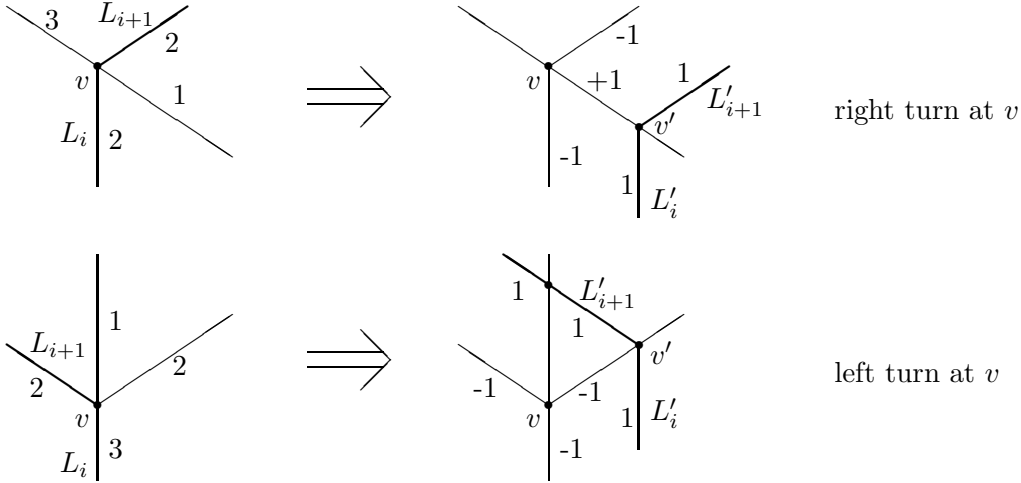
  \begin{picture}(145,65)
  \put(0,35){\begin{picture}(145,30)
\put(12,20){\circle*{1.0}}               
\put(12,20){\line(3,-2){18}}
\put(12,20){\line(-3,2){12}}
   \thicklines
\put(12,20){\line(0,-1){16}}
\put(12,20){\line(3,2){12}}
   \thinlines
\put(9,17){$v$}
\put(7,10){$L_i$}
\put(12,26){$L_{i+1}$}
\put(5,25){3}
\put(21,22){2}
\put(13.5,9){2}
\put(22,15){1}
\put(40,17){\line(1,0){10}}
\put(40,15){\line(1,0){10}}
\put(51,16){\line(-1,1){4}}
\put(51,16){\line(-1,-1){4}}
\put(72,20){\circle*{1.0}}               
\put(84,12){\circle*{1.0}}
\put(72,20){\line(0,-1){16}}
\put(72,20){\line(3,-2){18}}
\put(72,20){\line(-3,2){12}}
\put(72,20){\line(3,2){12}}
   \thicklines
\put(84,12){\line(0,-1){12}}
\put(84,12){\line(3,2){12}}
   \thinlines
\put(69,17){$v$}
\put(86,11){$v'$}
\put(85.5,1){$L'_i$}
\put(93,14.5){$L'_{i+1}$}
\put(73.5,7){-1}
\put(77,17){+1}
\put(81,23){-1}
\put(81.5,3){1}
\put(89,18){1}
\put(110,13){right turn at $v$}
  \end{picture}}
  \put(0,0){\begin{picture}(145,30)
\put(12,10){\circle*{1.0}}               
\put(12,10){\line(0,1){20}}
\put(12,10){\line(3,2){18}}
   \thicklines
\put(12,10){\line(0,-1){10}}
\put(12,10){\line(-3,2){12}}
   \thinlines
\put(9,7){$v$}
\put(7,1){$L_i$}
\put(3,17){$L_{i+1}$}
\put(3,11){2}
\put(22,14){2}
\put(13.5,3){3}
\put(13.5,20){1}
\put(40,17){\line(1,0){10}}
\put(40,15){\line(1,0){10}}
\put(51,16){\line(-1,1){4}}
\put(51,16){\line(-1,-1){4}}
\put(72,10){\circle*{1.0}}               
\put(84,18){\circle*{1.0}}
\put(72,26){\circle*{1.0}}
\put(72,10){\line(0,-1){10}}
\put(72,10){\line(0,1){20}}
\put(72,10){\line(-3,2){12}}
\put(72,10){\line(3,2){18}}
   \thicklines
\put(84,18){\line(0,-1){14}}
\put(84,18){\line(-3,2){18}}
   \thinlines
\put(69,7){$v$}
\put(85.5,15.5){$v'$}
\put(85.5,8){$L'_i$}
\put(76,24){$L'_{i+1}$}
\put(73.5,3){-1}
\put(62,11){-1}
\put(76,10.5){-1}
\put(81.5,7){1}
\put(67,24){1}
\put(75,19){1}
\put(110,13){left turn at $v$}
  \end{picture}}
  \end{picture}
 \end{center}

  \begin{lemma}  \label{lm:pp}
There exists $\bar\eps_1$, $0<\bar\eps_1\le\bar\eps_0$ such that
$\Pscr'$ is a pre-honeycomb for any nonnegative real
$\eps\le\bar\eps_1$.
  \end{lemma}
  \begin{proof}
Let $0<\eps<\bar\eps_0$. To see that $\Pscr'$ has zero tension
everywhere, it suffices to check this property at the bend vertices
$u_i$ of $P$ and their ``copies'' $u'_i$. For a bend vertex $u_i$,
let $\Pscr_i$ be the $\Xi$-system formed by the lines $L_i,L_{i+1},R_i$
with weights $-1,-1,w'(R_i)$, respectively, and $\Pscr'_i$ the
$\Xi$-system formed by the lines $L'_i,L'_{i+1},R_i$ with weights
$1,1,w'(R_i)$, respectively. One can see that
$\Pscr_i$ has zero tension at $u_i$ and $\Pscr'_i$ has zero tension
at $u'_i$, wherever (right or left) $P$ turns at $u_i$.
This implies the zero tension property for $\Pscr'$, taking into
account that $\Pscr$ has this property (as $\Pscr_i$ and $\Pscr'_i$
describe the corresponding local changes concerning $u_i,u'_i$ when
$\Pscr$ turns into $\Pscr'$).

It remains to explain that the numbers $(w')^s_j(v)$ are nonnegative
for all corresponding $v,j,s$ when $\eps>0$ is sufficiently small.

By (i),(ii) in Lemma~\ref{lm:lpath},
$\tilde w\ge 0$ for all $e\in\Escr$. So the only
case when $(w')^s_p(v)$ might be negative is when $v$ lies on a line
$R_i$ with weight $-1$ and $R_i$ is a $\Xi_p$-line.
Note that if $u_j=u_{j'}$ for some $j\ne j'$,
i.e., $P$ turns at the corresponding vertex $v$ of
$\Hscr$ twice, then the points $u'_j$ and $u'_{j'}$ move along
different rays out of $v$ (this can be concluded from~\refeq{shiftP},
taking into account (i) in Lemma~\ref{lm:lpath}). Hence we can choose
$\eps>0$ such that the interiors of the lines $R_i$ are pairwise
disjoint.

Consider $R_i$ with $w'(R_i)=-1$, and let $e$ be the edge dominating
at the vertex $v=u_i$ and different from $q_{t(i)}$ and $q_{t(i)+1}$.
Observe that the point $u'_i$ moves just along $e$,
so the line $R_i$ is entirely contained in $e$ when
$\eps$ does not exceed the length of $e$.
We show that $\tilde w(e)>0$, whence the result follows.
This is equivalent to saying that the number $\alpha$
of lines among $L_1,\ldots,L_r$ that contain $e$ (equal to 0,1 or 2) is
strictly less than $w(e)$. For a contradiction, suppose $\alpha=w(e)$
(as $\alpha\le w(e)$, by Lemma~\ref{lm:lpath}).
This implies that the number of bends of $P$ at $v$ using the edge $e$
(equal to the number of lines $L_j$ that contain $e$ and have one end
at $v$) is at least $|\diver_w(v)|$. Then the total number of bends
at $v$ is greater than $|\diver_w(v)|$ (as the bend
$(q_{t(i)},u_i,q_{t(i)+1})$ does not use $e$), contradicting (iii) in
Lemma~\ref{lm:lpath}. So $\alpha<w(e)$, as required.
   \end{proof}

For $\eps$ as in the lemma, $\Pscr'$ determines a
honeycomb $\Hscr'$, as explained in Section~\ref{sec:honey}.
We say that $\Hscr'$ is obtained from $\Hscr$ by the
{\em right $\eps$-deformation} of $P$.

\section{Proof of the Theorem} \label{sec:proof}

In what follows, given a honeycomb with edge weights $w''$, we say
that vertices $u,v$ have {\em the same sign} if
$\diver_{w''}(u)\diver_{w''}(v)\ge 0$, and call $|\diver_{w''}(v)|$
the {\em excess} at $v$.

Consider a concave cocirculation $h$ in a convex grid $G$ and the
honeycomb $\Hscr=(\Vscr,\Escr,w)$ determined by $h$. Define $\beta=
\beta_\Hscr$ to be the total weight of {\em nonintegral} semiinfinite
edges of $\Hscr$, $\delta=\delta_\Hscr$ to be the total excess of
{\em nonintegral} vertices of $\Hscr$, and $\omega=\omega_\Hscr$ to be
the total weight of edges incident to {\em integral} vertices.
We prove Theorem~\ref{tm:main} by induction on
   $$
\eta:=\eta_\Hscr:=\beta+\delta-\omega,
   $$
considering all concave cocirculations $h$ in the given $G$. Observe
that $\omega$ does not exceed the number of edges of $G$;
hence $\eta$ is bounded from below. Note also that in the case
$\beta=\delta=0$, all edges of $\Hscr$ are integral (whence $h$ is
integral as well). Indeed, suppose the set $\Ebullet$ of nonintegral
edges of $\Hscr$ is nonempty, and let $e\in\Ebullet$. Take the maximal
line $L$ that contains $e$ and is covered by edges of $\Hscr$.
Since $d^c(L)=d^c(e)$ is not an integer and $\beta=0$,
$L$ contains no semiinfinite edge; so $L$ is finite. The maximality
of $L$ implies that each end $v$ of $L$ is a vertex of $\Hscr$ and,
moreover, $\diver_w(v)\ne 0$. Also $v\in\Vbullet$. Then $\delta\ne 0$;
a contradiction.

Thus, we may assume that $\beta+\delta>0$ and that the theorem is
valid for all concave cocirculations on $G$ whose corresponding
honeycombs $\Hscr'$ satisfy $\eta_{\Hscr'}<\eta_\Hscr$.
We use notation, constructions and facts from
Sections~\ref{sec:l_path},\ref{sec:deform}.

Choose $P=(v_0,q_1,v_1,\ldots,q_k,v_k)$ as in Lemma~\ref{lm:lpath}.
Note that if $P$ is a cycle, then the fact that all bends of $P$ are
of the same degree $120^\circ$ implies that there are two
consecutive bend vertices $u_i,u_{i+1}$ where the direction of turn of
$P$ is the same.
We are going to apply to $P$ the right $\eps$-deformation, assuming
that either $P$ is an open path, or $P$ is a cycle having two
consecutive bend vertices where it turns right (for $P$ can be
considered up to reversing).

We gradually grow the parameter $\eps$ from zero,
obtaining the (parameteric) honeycomb $\Hscr'=(\Vscr',\Escr',w')$ as
described in Section~\ref{sec:deform}. Let $\bar\eps_1$ be specified
as the {\em maximum} real or $+\infty$, with $\bar\eps_1\le\bar\eps_0$,
satisfying the assertion of Lemma~\ref{lm:pp}.
(Such an $\bar\eps_1$ exists, as if $\eps'>0$ and
if $\Pscr'$ is a pre-honeycomb for any $0<\eps<\eps'$, then $\Pscr'$
is a pre-honeycomb for $\eps=\eps'$ either, by continuity and
compactness.) We stop growing $\eps$ as soon as it reaches the bound
$\bar\eps_1$ or at least one of the following events happens:

(E1) when $P$ is an open path, the constant coordinate of some of the
(semiinfinite or infinite) lines $L'_1$ and $L'_r$ becomes an integer;

(E2) two vertices of $\Hscr'$ with different signs meet (merge);

(E3) some line $L'_i$ meets an integer vertex $v$ of the original
$\Hscr$.

By the above assumption, if $P$ is a cycle, then $\bar\eps_1\le
\bar\eps_0<\infty$ (cf.~\refeq{length}).
And if $P$ is an open path, then the growth of
$\eps$ is bounded because of (E1). So we always finish with a finite
$\eps$; let $\bar\eps$ and $\bar\Hscr$ denote the resulting $\eps$ and
honeycomb, respectively. We assert that $\eta_{\bar\Hscr}<\eta_\Hscr$.
First of all notice that $\beta_{\bar\Hscr}\le\beta_\Hscr$ (in view of
(E1) and since the edges of $P$ are nonintegral).
Our further analysis relies on the following observations.

\smallskip
(i) When $\eps$ grows, each point $u'_i$ uniformly moves along a ray
from $u_i$, and each line $L'_i$ uniformly moves in a normal
direction to $L_i$. This implies that one can select a finite sequence
$0=\eps(0)<\eps(1)<\ldots<\eps(N)=\bar\eps$ such that $N=O(|\Vscr|^2)$,
and for $t=0,\ldots,N-1$, the
honeycomb $\Hscr'$ does not change topologically when the parameter
$\eps$ ranges over the open interval $(\eps(t),\eps(t+1))$.

\smallskip
(ii) When $\eps$ starts growing from zero, each vertex $v$ of $\Hscr$
occurring in $P$ splits into several vertices whose total divergency
is equal to the original divergency at $v$. By the construction of
$\Pscr'$ and (iii) in Lemma~\ref{lm:lpath}, these vertices have the
same sign. This implies that the total excess of these vertices is
equal to the original excess at $v$. Each arising vertex $u'_i$ has
excess 1, which preserves during the process except possibly for
those moments $\eps(t)$ when $u'_i$ can meet another vertex of
$\Hscr'$. When two or more vertices meet, their divergencies are
added up. Therefore, the sum of their excesses (before the meeting)
is strictly more than the resulting excess if some of these
vertices have different signs. This implies that $\delta$ reduces
if (E2) happens, or if $\eps$ reaches $\bar\eps_0$
(since the ends of any finite line $L'_i$ have different signs, and
the vertices $u'_{i-1}$ and $u'_i$ meet when $L'_i$ vanishes). The
latter situation is illustrated in the picture.
 \begin{center}
  \unitlength=1mm
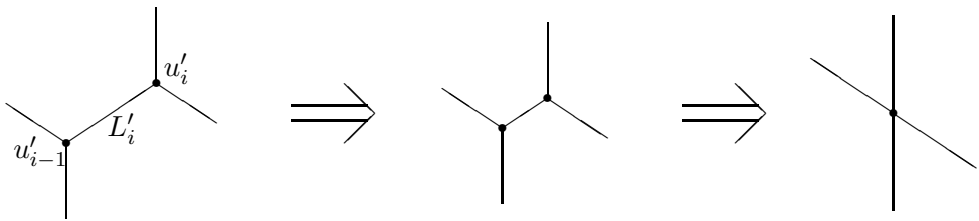
  \begin{picture}(140,28)
\put(0,0){\begin{picture}(51,28)
\put(10,10){\circle*{1.0}}               
\put(10,10){\line(-3,2){8}}
\put(10,10){\line(0,-1){10}}
\put(10,10){\line(3,2){12}}
\put(22,18){\circle*{1.0}}
\put(22,18){\line(3,-2){8}}
\put(22,18){\line(0,1){10}}
\put(3,8){$u'_{i-1}$}
\put(23,19){$u'_i$}
\put(15.5,11){$L'_i$}
\put(40,13){\line(1,0){10}}
\put(40,15){\line(1,0){10}}
\put(51,14){\line(-1,1){4}}
\put(51,14){\line(-1,-1){4}}
  \end{picture}}
\put(55,0){\begin{picture}(51,28)
\put(13,12){\circle*{1.0}}               
\put(13,12){\line(-3,2){8}}
\put(13,12){\line(0,-1){10}}
\put(13,12){\line(3,2){6}}
\put(19,16){\circle*{1.0}}
\put(19,16){\line(3,-2){8}}
\put(19,16){\line(0,1){10}}
\put(37,13){\line(1,0){10}}
\put(37,15){\line(1,0){10}}
\put(48,14){\line(-1,1){4}}
\put(48,14){\line(-1,-1){4}}
  \end{picture}}
\put(120,14){\circle*{1.0}}               
\put(120,14){\line(0,1){13}}
\put(120,14){\line(0,-1){13}}
\put(120,14){\line(3,-2){11}}
\put(120,14){\line(-3,2){11}}
  \end{picture}
 \end{center}

(iii) Let $v$ be an integral vertex of the initial honeycomb $\Hscr$
and let $W$ be the total weight of its incident edges in a current
honeycomb $\Hscr'$ (depending on $\eps$).
If, at some moment, the point $v$ is captured by the interior of
some line $L'_i$, this increases $W$ by 2. Now suppose some vertex
$u'_i$ meets $v$. If $P$ turns right at $u_i$, then $W$ increases
by at least $w'(R_i)=1$. And if $P$ turns left at $u_i$, then the
lines $L'_i,L'_{i+1}$ do not vanish (cf.~\refeq{length}), whence $W$
increases by $w'(L'_i)+w'(L'_{i+1})+w'(R_i)=1$. Therefore, $\omega$
increases when (E3) happens.

\smallskip
(iv) Let $\eps=\bar\eps_1<\bar\eps_0$. By the maximality of
$\bar\eps_1$ and
reasonings in the proof of Lemma~\ref{lm:pp}, a growth of $\eps$
beyond $\bar\eps_1$ would make some value $(w')^s_i(v)$ be
negative. This can happen only in two cases: (a) some line $R_j$ with
weight $-1$ is covered by edges of $\Hscr$ when $\eps=\bar\eps_1$,
and not covered when $\eps>\bar\eps_1$; or (b) some $R_j,R_{j'}$ with
weight $-1$ each lie on the same infinite line, the points
$u'_j$ and $u'_{j'}$ move toward each other when $\eps<\bar\eps_1$ and
these points meet when $\eps$ reaches $\bar\eps_1$ (then $R_j,R_{j'}$
become overlapping for $\eps>\bar\eps_1$). One can see that, in case
(a), $u'_j$ meets a vertex $v$ of $\Hscr$ (when $\eps=\bar\eps_1$)
and the signs of $v,u_j$ are different, and in case (b), the signs
of $u_j,u_{j'}$ are different as well. In both cases, $\delta$
decreases (cf.~(ii)). Case (b) is illustrated in the picture.
 \begin{center}
  \unitlength=1mm
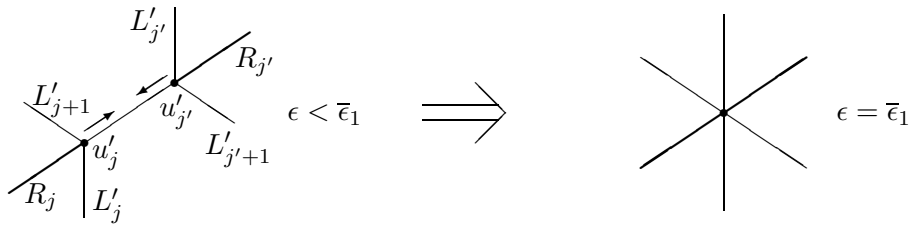
  \begin{picture}(120,28)
\put(10,10){\circle*{1.0}}               
\put(10,10){\line(-3,2){8}}
\put(10,10){\line(0,-1){10}}
\put(10,10){\line(3,2){12}}
\put(22,18){\circle*{1.0}}
\put(22,18){\line(3,-2){8}}
\put(22,18){\line(0,1){10}}
\put(10,11.5){\vector(3,2){4}}
\put(21,19){\vector(-3,-2){4}}
\put(11,8){$u'_j$}
\put(20,13.5){$u'_{j'}$}
\put(2,2){$R_j$}
\put(11,1){$L'_j$}
\put(3,15){$L'_{j+1}$}
\put(30,20){$R_{j'}$}
\put(16,25){$L'_{j'}$}
\put(26,8.5){$L'_{j'+1}$}
\put(37,13){$\eps<\bar\eps_1$}
\thicklines
\put(10,10){\line(-3,-2){10}}
\put(22,18){\line(3,2){10}}
\thinlines
\put(55,13){\line(1,0){10}}
\put(55,15){\line(1,0){10}}
\put(66,14){\line(-1,1){4}}
\put(66,14){\line(-1,-1){4}}
\put(95,14){\circle*{1.0}}               
\put(95,14){\line(0,1){13}}
\put(95,14){\line(0,-1){13}}
\put(95,14){\line(3,-2){11}}
\put(95,14){\line(-3,2){11}}
\put(110,13){$\eps=\bar\eps_1$}
\thicklines
\put(95,14){\line(-3,-2){11}}
\put(95,14){\line(3,2){11}}
  \end{picture}
 \end{center}

Using these observations, one can conclude that, during the process,
$\beta$ and $\delta$ are monotone nonincreasing and $\omega$ is
monotone nondecreasing. Moreover, at least one of these values
must change. Hence $\eta_{\bar\Hscr}<\eta_\Hscr$, as required.
(We leave it to the reader to examine details more carefully where
needed.)

Let $\bar h$ be the concave cocirculation in $G$ determined by
$\bar\Hscr$.
(The graph $G$ does not change as it is determined by the list of
numbers $w(\Bscr^s_i)$, defined in
Section~\ref{sec:honey}, and this list preserves.)
To justify the induction and finish the proof of the theorem,
it remains to explain that
   \begin{equation} \label{eq:OI}
I_h\subseteq I_{\bar h} \qquad \mbox{and}
\qquad O_h\subseteq O_{\bar h}.
   \end{equation}

Denote by $\Hscr_\eps$ and $h_\eps$ the current honeycomb $\Hscr'$ and
the induced concave cocirculation $h'$ at a moment $\eps$ in the
process, respectively.
The correspondence between the vertices of $\Hscr_\eps$ and the
flatspaces of $h_\eps$ (explained in Section~\ref{sec:honey}) implies
that for each edge $e\in E(G)$, the function $h_\eps(e)$ is continuous
within each interval $(\eps(t),\eps(t+1))$ (cf. (i) in the
above analysis). We assert that $h'=h_\eps$ is continuous in the
entire segment $[0,\bar\eps]$ as well.

To see the latter, consider the honeycomb $\Hscr_{\eps(t)}$ for
$0\le t<N$. When $\eps$ starts growing from $\eps(t)$
(i.e., $\eps(t)<\eps<\eps(t+1)$ and $\eps-\eps(t)$ is small)
the set $Q(v)$ of vertices of $\Hscr'=\Hscr_\eps$ arising from a
vertex $v$ of $\Hscr_{\eps(t)}$ (by splitting or moving or preserving
$v$) is located in a small neighbourhood of $v$.
Moreover, for two distinct vertices $u,v$ of $\Hscr_{\eps(t)}$,
the total weight of edges of $\Hscr_\eps$ connecting $Q(u)$ and
$Q(v)$ is equal to the weight of the edge between $u$ and $v$ in
$\Hscr_{\eps(t)}$ (which is zero when the edge does not exist),
and all these edges are parallel. This implies that
for each vertex $v$ of $\Hscr_{\eps(t)}$, the arising subgrids
$G_{v'}$, $v'\in Q(v)$, in $G$ (concerning $h_\eps$) give a partition
of the subgrid $G_v$ (concerning $h_{\eps(t)}$), i.e., the set of
little triangles of $G$ contained in these $G_{v'}$ coincides
with the set of little triangles occurring in $G_v$. So $h'$ is
continuous within $[\eps(t),\eps(t+1))$.

Similarly, for each vertex $v$ of $\Hscr_{\eps(t)}$ ($0<t\le N$),
the subgrid $G_v$ is obtained by gluing together the subgrids $G_{v'}$,
$v'\in Q'(v)$, where $Q'(v')$ is the set of vertices of $\Hscr_\eps$
(with $\eps(t-1)<\eps<\eps(t)$) which produce $v$ when $\eps$ tends to
$\eps(t)$. So $h'$ is continuous globally.
Now since no integral vertex of the initial honeycomb can move or
split during the process (but merely merge with another vertex of
$\Hscr'$ if (E3) happens), we conclude that the cocirculation
preserves in all little triangles where it is integral initially,
yielding the first inclusion in~\refeq{OI}.

The second inclusion in~\refeq{OI} is shown in a similar fashion,
relying on (E1).

This completes the proof of Theorem~\ref{tm:main}.

\section{Polyhedra $\Cscr(G,h_0)$ Having Vertices with Big
Denominators} \label{sec:frac}

In this section we describe a construction proving
Theorem~\ref{tm:frac}. We start with some definitions and auxiliary
statements.

Given a concave cocirculation $h$ in a concave grid $G$, the {\em
tiling} $\tau_h$ is the subdivision of the polygon $\Rscr(G)$ spanned
by $G$ into the flatspaces $T$ of $h$. We also call $T$ a {\em tile}
in $\tau_h$. The following elementary property of tilings will be
important for us.

  \begin{lemma}  \label{lm:tiling}
Let $F\subseteq E(G)$, $h_0:F\to\Rset$, and $h\in\Cscr(G,h_0)=:
\Cscr$. Then $h$ is a vertex of $\Cscr$ if and only if $h$ is uniquely
determined by its tiling, i.e., $h'\in\Cscr$ and $\tau_{h'}=\tau_h$
imply $h'=h$.
  \end{lemma}
\begin{proof}
For a little rhombus $\rho$ of $G$, the fact that the sum of values of
$h$ over each of the two 3-circuits in $\rho$ is zero implies that
if rhombus inequality~\refeq{concave} turns into equality for one pair
of parallel edges in $\rho$, then it does so for the other pair. Since
$\Cscr$ is described by a system of linear constraints where the
inequalities are exactly of the form~\refeq{concave}, $h$ is a vertex
of $\Cscr$ if and only if it is uniquely determined by the set $Q$ of
little rhombi for which~\refeq{concave} turns into equality.
Now observe that each tile in $\tau_h$ one-to-one corresponds to a
component of the graph whose vertices are the little triangles of $G$
and whose edges are the pairs of little triangles forming rhombi in
$Q$. This implies the lemma.
  \end{proof}

We will use a re-formulation of this lemma involving honeycombs.
Let us say that honeycombs $\Hscr=(\Vscr,\Escr,w)$ and
$\Hscr'=(\Vscr',\Escr',w')$ are {\em conformable} if
$|\Vscr|=|\Vscr'|$, $|\Escr|=|\Escr'|$, and there are bijections
$\alpha:\Vscr\to\Vscr'$ and $\beta:\Escr\to\Escr'$ such that: for each
vertex $v\in\Vscr$ and each edge $e\in\Escr$ incident to $v$, the edge
$\beta(e)$ is incident to $\alpha(v)$, $w(e)=w'(\beta(e))$, and if $e$
is contained in $\Xi^s_i(v)$, then $\beta(e)$ is contained in
$\Xi^s_i(\alpha(v))$. (This matches the situation when two concave
cocirculations in $G$ have the same tiling.)

Next, for $\Fscr\subseteq\Escr$, we call $\Hscr$ {\em extreme} with
respect to $\Fscr$, or $\Fscr$-{\em extreme}, if there is no honeycomb
$\Hscr'\ne\Hscr$ such that $\Hscr'$ is conformable to $\Hscr$ and
satisfies $d^c(\beta(e))=d^c(e)$ for all $e\in\Fscr$, where $\beta$
is the corresponding bijection.
  \Xcomment{
(This is analogous to the property that a concave cocirculation is a
vertex of the corresponding polyhedron. Note that in contrast to the
concave cocirculation case where $F$ is a subset of edges of a fixed
grid $G$, the above definition of extreme honeycomb involves a subset
of edges of the honeycomb itself. Using some formalism, one can give a
definition not dependent on a particular honeycomb (e.g., in spirit
of~\cite{KT}) but we do not need it in this paper.)
  }
Then the relationship between the tiling of concave cocirculations and
the vertex sets of honeycombs leads to the following re-formulation of
Lemma~\ref{lm:tiling}.

\begin{corollary}  \label{cor:extreme}
Let $F\subseteq E(G)$, $h_0:F\to\Rset$, and $h\in\Cscr(G,h_0)$. Let
$\Hscr$ be the honeycomb determined by $h$ and let $\Fscr$ be the
subset of edges of $\Hscr$ corresponding to sides of tiles in $\tau_h$
that contain at least one edge from $F$. Then $h$ is a vertex of
$C(G,h_0)$ if and only if $\Hscr$ is $\Fscr$-extreme.
  \end{corollary}

One sufficient condition on extreme honeycombs will be used later.
Let us say that a line $\ell$ in $\Rset^2$ is a {\em line of} $\Hscr$
if $\ell$ is covered by edges of $\Hscr$. Then
  \begin{itemize}
\item[(C)]
$\Hscr$ is $\Fscr$-extreme if each vertex of $\Hscr$ is contained in
at least two maximal lines of $\Hscr$, each containing an edge from
$\Fscr$.
  \end{itemize}
Indeed, if two different maximal lines $L,L'$ of $\Hscr$ intersect
at a vertex $v$, then the constant coordinates of $L,L'$
determine the dual coordinates of $v$. One can see that if $\Hscr'$ is
conformable to $\Hscr$, then the images of $L,L'$ in $\Hscr'$ are
maximal lines there and they intersect at the image of $v$. This
easily implies~(C).

The idea of our construction is as follows. Given a positive
integer $k$, we will devise two honeycombs. The first honeycomb
$\Hscr'=(\Vscr',\Escr',w')$ has the following properties:
 \begin{itemize}
\item[(P1)]
 \begin{itemize}
\item[(i)] the boundary $\Bscr(\Hscr')$ is partitioned into three sets
$\Bscr'_1,\Bscr'_2,\Bscr'_3$, where $\Bscr'_i$ consists of the
semiinfinite edges of the form $\Xi^+_i(\cdot)$,
and $w'(\Bscr'_i)\le Ck$, where $C$ is a constant;
\item[(ii)] the constant coordinates of all edges of $\Hscr'$ are
integral;
\item[(iii)] $\Hscr'$ is extreme with respect to $\Bscr'_1\cup
\Bscr'_2$.
  \end{itemize}
  \end{itemize}

The second honeycomb $\Hscr''=(\Vscr'',\Escr'',w'')$ has the following
properties:
 \begin{itemize}
\item[(P2)]
 \begin{itemize}
\item[(i)] each semiinfinite edge of $\Hscr''$ is contained in a line
of $\Hscr'$ (in particular, $d^c(e)$ is an integer for each
$e\in\Bscr(\Hscr'')$), and $w''(\Bscr(\Hscr''))\le w'(\Bscr(\Hscr'))$;
\item[(ii)] there is $e\in\Escr''$ such that the denominator of
$d^c(e)$ is equal to $k$;
\item[(iii)] $\Hscr'$ is extreme with respect to its boundary
$\Bscr(\Hscr'')$.
  \end{itemize}
  \end{itemize}

Define the {\em sum} $\Hscr'+\Hscr''$ to be the
honeycomb $\Hscr=(\Vscr,\Escr,w)$ determined by the pre-honeycomb
in which the line set is the (disjoint) union of $\Escr'$ and $\Escr''$,
and the weight of a line $e$ is equal to $w'(e)$ for $e\in\Escr'$, and
$w''(e)$ for $e\in\Escr''$. Then each edge of $\Hscr$ is contained in
an edge of $\Hscr'$ or $\Hscr''$, and conversely, each edge of
$\Hscr'$ or $\Hscr''$ is contained in a line of $\Hscr$. Using this,
one can derive from (P1) and (P2) that
 \begin{itemize}
\item[(P3)]
 \begin{itemize}
\item[(i)] the boundary $\Bscr(\Hscr)$ is partitioned into three sets
$\Bscr_1,\Bscr_2,\Bscr_3$, where $\Bscr_i$ consists of the
semiinfinite edges of the form $\Xi^+_i(\cdot)$,
and $w(\Bscr_i)\le 2Ck$;
\item[(ii)] $d^c(e)\in\Zset$ for all $e\in\Bscr(\Hscr)$, and $d^c(e)$
has denominator $k$ for some edge $e$ of $\Hscr$;
\item[(iii)] $\Hscr$ is extreme with respect to $\Bscr_1\cup\Bscr_2$.
  \end{itemize}
  \end{itemize}
(Property (iii) follows from (P1)(iii) and (P2)(i),(iii).)

Now consider the grid $G$ and the concave cocirculation $h$ in it
determined by $\Hscr$. By (P3)(i), $G$ is a 3-side grid of size
at most $2Ck$, with sides $B_1,B_2,B_3$ corresponding to
$\Bscr_1,\Bscr_2,\Bscr_3$, respectively. Let $F:=B_1\cup B_2$ and let
$h_0$ be the restriction of $h$ to $F$ (considering a side as edge
set). Then (P3) together with Corollary~\ref{cor:extreme} implies that
$G,h_0,h$ are as required in Theorem~\ref{tm:frac}.

\medskip
It remains to devise $\Hscr'$ and $\Hscr''$ as above. To devise
$\Hscr'$ is rather easy. It can be produced by truncating the dual
grid formed by all lines with integer constant coordinates. More
precisely, let $n$ be a positive integer (it will depend on $k$). For
a point $x\in\Rset^2$ and $i=1,2,3$, let $x^i$ stand for the dual
coordinate $d_i(x)$. The vertex set $\Vscr'$ consists of the points
$v$ such that
 \begin{gather*}
 v^i\in\Zset \quad\mbox{and}\quad |v^i|<n, \qquad i=1,2,3,\\
 v^1-v^2\le n,\quad  v^2-v^3\le n,\quad  v^3-v^1\le n.
 \end{gather*}
The finite edges of $\Hscr'$ have unit weights and connect the pairs
$u,v\in\Vscr'$ such that
  $$
 |u^1-v^1|+|u^2-v^2|+|u^3-v^3|=2.
  $$
The semiinfinite edges and their weights are assigned as follows
(taking indices modulo 3):
 \begin{itemize}
\item
if $v\in\Vscr'$, $i\in\{1,2,3\}$, and $v^i-v^{i+1}\in\{n,n-1\}$, then
$\Hscr'$ has edge $e=\Xi^+_{i-1}(v)$, and the weight of $e$ is equal
to 1 if $v^i-v^{i+1}=n-1$ and $v^i,v^{i+1}\ne 0$, and equal to 2
otherwise.
  \end{itemize}
The case $n=3$ is illustrated in the picture, where the semiinfinite
edges with weight 2 are drawn in bold.
 \begin{center}
  \unitlength=1mm
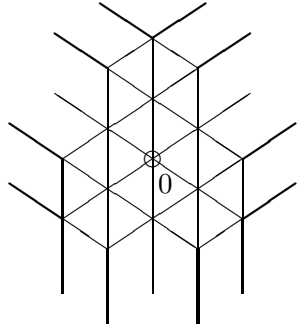
  \begin{picture}(45,45)
\put(22,22){\circle{2.0}}
\put(10,14){\line(0,1){8}}
\put(16,10){\line(0,1){24}}
\put(22,38){\line(0,-1){34}}
\put(28,10){\line(0,1){24}}
\put(34,14){\line(0,1){8}}
\put(10,14){\line(3,2){25}}
\put(10,22){\line(3,2){18}}
\put(16,10){\line(3,2){18}}
\put(16,34){\line(3,2){6}}
\put(28,10){\line(3,2){6}}
\put(10,14){\line(3,-2){6}}
\put(10,22){\line(3,-2){18}}
\put(34,14){\line(-3,2){25}}
\put(16,34){\line(3,-2){18}}
\put(22,38){\line(3,-2){6}}
\put(22.7,17.5){0}
\thicklines
\put(10,14){\line(0,-1){10}}
\put(16,10){\line(0,-1){10}}
\put(28,10){\line(0,-1){10}}
\put(34,14){\line(0,-1){10}}
\put(10,14){\line(-3,2){7}}
\put(10,22){\line(-3,2){7}}
\put(16,34){\line(-3,2){7}}
\put(22,38){\line(-3,2){7}}
\put(22,38){\line(3,2){7}}
\put(28,34){\line(3,2){7}}
\put(34,22){\line(3,2){7}}
\put(34,14){\line(3,2){7}}
  \end{picture}
 \end{center}

One can check that $\Hscr'$ is indeed a honeycomb and satisfies
(P1)(i),(ii) (when $n=O(k)$).
Also each vertex belongs to two lines of $\Hscr'$, one
containing a semiinfinite edge in $\Bscr'_1$, and the other in
$\Bscr'_2$. Then $\Hscr'$ is $(\Bscr'_1\cup\Bscr'_2)$-expreme, by
assertion (C). So $\Hscr'$ satisfies (P1)(iii) as well. Note that
the set of semiinfinite edges of $\Hscr'$ is dense, in the sense that
for each $i=1,2,3$ and $d=-n+1,\ldots,n-1$, there is a boundary edge
$e$ of the form $\Xi^+_i(\cdot)$ such that $d^c(e)=d$.

\medskip
Next we have to devise $\Hscr''$ satisfying (P2), which is less
trivial. In order to facilitate the description and technical details,
we go in reverse direction: we construct a certain concave
cocirculation $\tilde h$ is a convex grid $\tilde G$ and then
transform it into the desired honeycomb.

The grid $\tilde G$ spans a hexagon with S- and N-sides of length 1
and with SW-, NW-, SE-, and NE-sides of length $k$. We denote the
vertices in the big sides (in the order as they follow in the
side-path) by:
 \begin{itemize}
  \item[v1.] $x_k,x_{k-1},\ldots,x_0\quad$ for the SW-side;
  \item[v2.] $x'_k,x'_{k-1},\ldots,x'_0\quad$ for the NW-side;
  \item[v3.] $y_0,y_1,\ldots,y_k\quad$ for the SE-side;
  \item[v4.] $y'_0,y'_1,\ldots,y'_k\quad$ for the NE-side.
 \end{itemize}
(Then $x_0=x'_0$ and $y_0=y'_0$.) We also distinguish the vertices
$z_i:=x_i+\xi_1$ and $z'_i:=x'_i+\xi_1$ for $i=1,\ldots,k$
(then $z_0=z'_0$, $z_k=y_k$, $z'_k=y'_k$.) We arrange an
(abstract) tiling $\tau$ in $\tilde G$. It is symmetric w.r.t. the
horizontal line $x_0y_0$ and consists of
  \begin{itemize}
\item[t1.]
$2k-2$ trapezoids and two little triangles obtained by subdividing the
rhombus $R:=y_ky_0y'_kz_0$ (labeled via its vertices) by the
horizontal lines passing $y_{k-1},\ldots, y_0,y'_1,\ldots,y'_{k-1}$;
  \item[t2.] the little rhombus $\bar\rho:=x_1z_0x'_1x_0$;
  \item[t3.]
$4k-2$ little triangles $\Delta_i:=x_iz_iz_{i-1}$, $\nabla'_i:=
x'_iz'_{i-1}z'_i$ for $i=1,\ldots,k$, and $\nabla_j:=x_jx_{j+1}z_j$,
$\Delta'_j:=x'_jz'_jx'_{j+1}$ for $j=1,\ldots,k-1$.
  \end{itemize}

Define the function $h_0:B(\tilde G)\to\Zset$ by:
  \begin{itemize}
 \item[h1.] $h_0(x_ix_{i-1}):=h_0(x'_ix'_{i-1}):=i-1$ for
            $i=2,\ldots,k$;
 \item[h2.] $h_0(x_1x_0):=h_0(x'_1x'_0):=-1$ and
           $h_0(x_ky_k):=h_0(x'_ky'_k):=0$;
 \item[h3.] $h_0(y_{i-1}y_i):=h_0(y'_{i-1}y'_i):=-i+1$ for
            $i=1,\ldots,k$.
  \end{itemize}
The constructed $\tau$ and $h_0$ for $k=3$ are illustrated in the
picture (the numbers of the boundary edges indicate the values of
$h_0$).
 \begin{center}
  \unitlength=1mm
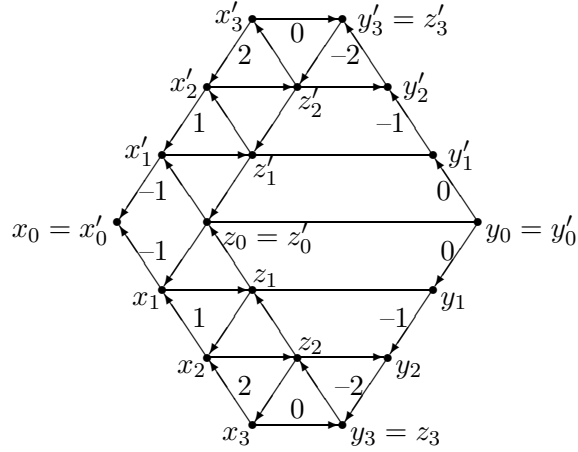
  \begin{picture}(80,60)
\put(38,3){\circle*{1.0}}
\put(50,3){\circle*{1.0}}
\put(32,12){\circle*{1.0}}
\put(44,12){\circle*{1.0}}
\put(56,12){\circle*{1.0}}
\put(26,21){\circle*{1.0}}
\put(38,21){\circle*{1.0}}
\put(62,21){\circle*{1.0}}
\put(20,30){\circle*{1.0}}
\put(32,30){\circle*{1.0}}
\put(68,30){\circle*{1.0}}
\put(26,39){\circle*{1.0}}
\put(38,39){\circle*{1.0}}
\put(62,39){\circle*{1.0}}
\put(32,48){\circle*{1.0}}
\put(44,48){\circle*{1.0}}
\put(56,48){\circle*{1.0}}
\put(38,57){\circle*{1.0}}
\put(50,57){\circle*{1.0}}
\put(38,3){\vector(1,0){11.5}}
\put(32,12){\vector(1,0){11.5}}
\put(26,21){\vector(1,0){11.5}}
\put(26,39){\vector(1,0){11.5}}
\put(32,48){\vector(1,0){11.5}}
\put(38,57){\vector(1,0){11.5}}
\put(44,12){\vector(1,0){11.5}}
\put(44,48){\vector(1,0){11.5}}
\put(38,3){\vector(-2,3){5.7}}
\put(50,3){\vector(-2,3){5.7}}
\put(32,12){\vector(-2,3){5.7}}
\put(44,12){\vector(-2,3){5.7}}
\put(26,21){\vector(-2,3){5.7}}
\put(38,21){\vector(-2,3){5.7}}
\put(32,30){\vector(-2,3){5.7}}
\put(38,39){\vector(-2,3){5.7}}
\put(44,48){\vector(-2,3){5.7}}
\put(68,30){\vector(-2,3){5.7}}
\put(62,39){\vector(-2,3){5.7}}
\put(56,48){\vector(-2,3){5.7}}
\put(44,12){\vector(-2,-3){5.7}}
\put(38,21){\vector(-2,-3){5.7}}
\put(32,30){\vector(-2,-3){5.7}}
\put(26,39){\vector(-2,-3){5.7}}
\put(38,39){\vector(-2,-3){5.7}}
\put(32,48){\vector(-2,-3){5.7}}
\put(44,48){\vector(-2,-3){5.7}}
\put(38,57){\vector(-2,-3){5.7}}
\put(50,57){\vector(-2,-3){5.7}}
\put(56,12){\vector(-2,-3){5.7}}
\put(62,21){\vector(-2,-3){5.7}}
\put(68,30){\vector(-2,-3){5.7}}
\put(38,21){\line(1,0){24}}
\put(32,30){\line(1,0){36}}
\put(38,39){\line(1,0){24}}
\put(34,1){$x_3$}
\put(28,10){$x_2$}
\put(22,19){$x_1$}
\put(6,28){$x_0=x'_0$}
\put(21,38.5){$x'_1$}
\put(27,47.5){$x'_2$}
\put(33,56.5){$x'_3$}
\put(51,1){$y_3=z_3$}
\put(57,10){$y_2$}
\put(63,19){$y_1$}
\put(69,28){$y_0=y'_0$}
\put(64,38){$y'_1$}
\put(58,47){$y'_2$}
\put(52,56){$y'_3=z'_3$}
\put(44,13){$z_2$}
\put(38,22){$z_1$}
\put(34,27){$z_0=z'_0$}
\put(38,36){$z'_1$}
\put(44,45){$z'_2$}
\put(43,4){0}
\put(36,7){2}
\put(30,16){1}
\put(23,25){--1}
\put(23,33){--1}
\put(30,42){1}
\put(36,51){2}
\put(43,54){0}
\put(49,7){--2}
\put(55,16){--1}
\put(63,25){0}
\put(62.5,33){0}
\put(54.5,42){--1}
\put(48.5,51){--2}
  \end{picture}
 \end{center}

The tiling $\tau$ has the property that $h_0$ (as well as any function
on $B(\tilde G)$) is extendable to at most one cocirculation $\tilde h$
in $\tilde G$ such that $\tilde h$ is {\em flat} within each tile $T$
in $\tau$, i.e., satisfies $\tilde h(e)=\tilde h(e')$ for any parallel
edges $e,e'$ in $T$. This is because we can compute $\tilde h$,
step by step, using two observations: (i) the values of $\tilde h$ on
two nonparallel edges in a tile $T$ determine $\tilde h$ on all edges
in $T$, and (ii) if a side $S$ of a polygon $P$ in $\tilde G$
(spanned by a subgraph of $\tilde G$) is entirely contained in a tile,
then the values of $\tilde h$ on the edges of $S$ are
detemined by the values on the other sides of $P$ (since $\tilde h$ is
constant on $S$ and $\tilde h(Q^c)=\tilde h(Q^a)$,
where $Q^c$ ($Q^a$) is the set of boundary edges of $P$ oriented
clockwise (resp. anticlockwise) around $P$).

We assign $\tilde h$ as follows. First assign $\tilde h(e):=h_0(e)$
for all $e\in B(\tilde G)$. Second assign $\tilde h(z_0x_1)
:=\tilde h(z_0x'_1):=-1$ (by (i) for the rhombus $\bar\rho$ in $\tau$).
Third, for an edge $e$ in the horizontal line $z_0y_0$, assign
 \begin{multline*}
 \tilde h(e):=\frac1k \left(\tilde h(z_0x_1)-\tilde h(x_kx_1)
+\tilde h(x_ky_k)-\tilde h(y_0y_k) \right)\\
=\frac1k \left(\vphantom{\tilde h} -1-(1+\ldots+(k-1))+0
      -(0-1-\ldots-(k-1)) \right) =\frac1k(-1)=-1/k
 \end{multline*}
(by (ii) for the pentagon $z_0x_1x_ky_ky_0$, taking into account that
the side $z_0y_0$ contains $k$ edges), where $\tilde h(uv)$ is the sum
of values of $\tilde h$ on the edges of a side $uv$. Therefore,
$\tilde h(e)=-1/k$ for all horizontal edges in the big rhombus $R$,
and we then can assign
$\tilde h(z_{i+1}z_i):=\tilde h(z'_{i+1}z'_i):=i+1/k$
for $i=0,\ldots,k-1$ (by applying (i) to the trapezoids and little
triangles of $\tau$ in $R$).
Fourth, repeatedly applying (i) to the little triangles
$\Delta_1,\nabla_1,\Delta_2,\ldots,\nabla_{k-1},\Delta_k$ (in this
order) that form the trapezoid $x_1x_ky_kz_0$ in which $\tilde h$ is
already known on all side edges, we determine $\tilde h$ on all
edges in this trapezoid; and similarly for the symmetric trapezoid
$z'_0y'_kx'_kx'_1$.

One can check that the obtained cocirculation $\tilde h$ is
well-defined, and moreover, it is a concave cocirculation with tiling
$\tau$ (a careful examination of the corresponding rhombus
inequalities is left to the reader). Since $\tilde h$ is computed
uniquely in the process, it is a vertex of
$\Cscr(\tilde G,h_0)$. Also $\tilde h$ is integral on the boundary
of $\tilde G$, has an entry with denominator $k$, and satisfies
$|\tilde h(e)|<2k$ for all $e\in E(\tilde G)$. The picture
indicates the values of $\tilde h$ in the horizontal stripe between
the lines $x_iy_i$ and $x_{i+1}y_{i+1}$ for $1\le i\le k-2$, and in
the horizontal stripe between $x_1y_1$ and $x'_1y'_1$.
 \begin{center}
  \unitlength=1mm
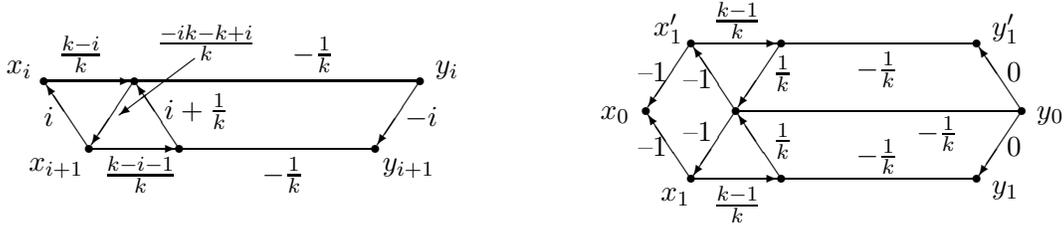
  \begin{picture}(140,35)
\put(12,14){\circle*{1.0}}               
\put(24,14){\circle*{1.0}}
\put(50,14){\circle*{1.0}}
\put(6,23){\circle*{1.0}}
\put(18,23){\circle*{1.0}}
\put(56,23){\circle*{1.0}}
\put(12,14){\vector(1,0){11.5}}
\put(6,23){\vector(1,0){11.5}}
\put(12,14){\vector(-2,3){5.7}}
\put(24,14){\vector(-2,3){5.7}}
\put(18,23){\vector(-2,-3){5.7}}
\put(56,23){\vector(-2,-3){5.7}}
\put(26,26){\vector(-4,-3){10}}
\put(24,14){\line(1,0){26}}
\put(18,23){\line(1,0){38}}
\put(4,11){$x_{i+1}$}
\put(51,11){$y_{i+1}$}
\put(1,24){$x_i$}
\put(58,24){$y_i$}
\put(14,9.5){$\frac{k-i-1}k$}
\put(22,18){$i+\frac1k$}
\put(35,9.5){$-\frac1k$}
\put(8,25.5){$\frac{k-i}k$}
\put(39,25.5){$-\frac1k$}
\put(21,27){$\frac{-ik-k+i}k$}
\put(6,17){$i$}
\put(54,17){$-i$}
\put(80,0){\begin{picture}(60,35)
\put(12,10){\circle*{1.0}}
\put(24,10){\circle*{1.0}}
\put(50,10){\circle*{1.0}}
\put(6,19){\circle*{1.0}}
\put(18,19){\circle*{1.0}}
\put(56,19){\circle*{1.0}}
\put(12,28){\circle*{1.0}}
\put(24,28){\circle*{1.0}}
\put(50,28){\circle*{1.0}}
\put(12,10){\vector(1,0){11.5}}
\put(12,28){\vector(1,0){11.5}}
\put(12,10){\vector(-2,3){5.7}}
\put(24,10){\vector(-2,3){5.7}}
\put(18,19){\vector(-2,3){5.7}}
\put(56,19){\vector(-2,3){5.7}}
\put(12,28){\vector(-2,-3){5.7}}
\put(18,19){\vector(-2,-3){5.7}}
\put(24,28){\vector(-2,-3){5.7}}
\put(56,19){\vector(-2,-3){5.7}}
\put(24,10){\line(1,0){26}}
\put(24,28){\line(1,0){26}}
\put(18,19){\line(1,0){38}}
\put(8,7){$x_1$}
\put(0,18){$x_0$}
\put(7,29){$x'_1$}
\put(52,8){$y_1$}
\put(58,18){$y_0$}
\put(52,29){$y'_1$}
\put(15,5.5){$\frac{k-1}k$}
\put(15,30){$\frac{k-1}k$}
\put(23,14){$\frac1k$}
\put(23,23){$\frac1k$}
\put(34,12){$-\frac1k$}
\put(34,23.5){$-\frac1k$}
\put(42,15){$-\frac1k$}
\put(5,13){--1}
\put(5,23){--1}
\put(11,15){--1}
\put(11,22){--1}
\put(54,13){0}
\put(54,23){0}
\end{picture}}
  \end{picture}
 \end{center}

Let $\tilde\Hscr$ be the honeycomb determined by
$(\tilde G,\tilde h)$. It is $\Bscr(\tilde\Hscr)$-extreme,
by Corollary~\ref{cor:extreme}, has all boundary edges integral and
has a finite edge with constant coordinate $1/k$.
We slightly modify $\tilde\Hscr$ in order to get rid of the
semiinfinite edges $e$ of the form $\Xi^-_i(v)$ in it. This is done
by truncating such an $e$ to finite edge $e'=vu$ and adding two
semiinfinite edges $a:=\Xi^+_{i-1}(u)$ and $b:=\Xi^+_{i+1}(u)$,
all with the weight equal to that of $e$ (which is 1 in our case).
Here $u$ is the integer point in $\Xi^-_i(v)\setminus\{v\}$ closest
to $v$. (Strictly
speaking, when applying this operation simultaneously to all such
edges $e$, we should handle the corresponding pre-honeycomb, as some
added edge may intersect another edge at an interior point.) The
obtained honeycomb $\Hscr''=(\Vscr'',\Escr'',w'')$ has all
semiinfinite edges $e$ of the form $\Xi^+_i(\cdot)$, with $d^c(e)$
being an integer between $-2k$ and $2k$. Also $w''(\Bscr(\Hscr''))<
2|B(\tilde G)|$. This honeycomb is extreme with respect to its
boundary since so is $\tilde\Hscr$ and since for $e',a,b$ as above,
the constant coordinate of $e'$ is determined by the constant
coordinates of $a,b$. Thus, $\Hscr''$ satisfies (P2) when $n>2k$
(to ensure that each semiinfinite edge of $\Hscr''$ is contained in
a line of $\Hscr'$).

Now $(G,h)$ determined by the honeycomb $\Hscr'+\Hscr''$ is as
required in Theorem~\ref{tm:frac}.

\medskip
\noindent
{\bf Remark 2.} In our construction of vertex $\tilde h$ of
$\Cscr(\tilde G,h_0)$, the design of tiling $\tau$ is borrowed from
one fragment in a construction in~\cite{DM} where one shows that
(in our terms) the polytope of semiconcave cocirculations in a 3-side
grid with fixed integer values on the boundary can have a vertex
with denominator $k$. Here we call a cocirculation
$h$ {\em semiconcave} if the inequality~\refeq{concave} holds
for each rhombus $\rho$ whose diagonal edge $d$ is parallel to $\xi_2$
or $\xi_3$ (but may be violated if $d$ is parallel to $\xi_1$).

\medskip
\noindent
{\bf Remark 3.} The first honeycomb $\Hscr'$ in our construction turns
out to be $\Fscr'$-extreme for many proper subsets $\Fscr'$ of
$\Bscr(\Hscr')$ and even of $\Bscr'_1\cup \Bscr'_2$.
For example, one can take $\Fscr':=\Bscr'_1\cup\{e\}$, where
$e$ is an arbitrary edge in $\Bscr'_2$ (a check-up is left to the
reader). Moreover, it is not difficult to produce triples of
boundary edges that possess such a property. For each of these sets
$\Fscr'$, the honeycomb $\Hscr=\Hscr'+\Hscr''$ is $\Fscr$-extreme,
where $\Fscr$ consists of the semiinfinite edges of $\Hscr$ contained
in members of $\Fscr'$. Then, by Corollary~\ref{cor:extreme},
the constructed concave cocirculation $h$ is a vertex of the
polyhedron $\Cscr(G,h\rest{F})$ as well, where $F$ is a subset of
boundary edges of $G$ whose images in $\Hscr$ cover $\Fscr$ (i.e., for
each edge $\ell\in\Fscr$ of the form $\Xi^+_i(\cdot)$, there is an
$e\in F\cap B_i$ with $h(e)=d^c(\ell)$). This gives a strengthening of
Theorem~\ref{tm:frac}.

\medskip
\noindent
{\bf Remark 4.}
One can try to describe the above construction implying
Theorem~\ref{tm:frac} so as to use the language of concave
cocirculations everywhere. However, this seems to be a more intricate
way. In particular, the operation on a pair $h',h''$ of concave
cocirculations analogous to taking the sum of honeycombs
$\Hscr',\Hscr''$ is less transparent (it is related to taking the
convolution of concave functions behind $h',h''$). This is why we
prefer to argue in terms of honeycombs.

\section{Concluding Remarks} \label{sec:concl}

The proof of Theorem~\ref{tm:main} in Section~\ref{sec:proof}
provides a strongly
polynomial algorithm which, given $G$ and $h$, finds $h'$ as required
in this theorem. Indeed, each parameter $\beta,\delta,\omega$ is
bounded by the number of edges of $G$, so the number of iterations
(viz. applications of the induction) is $O(n)$, where $n:=|E(G)|$.
As was explained, the number of moments $\eps$ when some line
$L'_i$ captures a vertex of $\Hscr$ or when two vertices $u'_i,u'_j$
meet is $O(|\Vscr|^2)$,
or $O(n^2)$, and these moments can be computed easily. To find
$\bar\eps_1$ is easy as well. Hence an iteration is performed in time
polynomial in $n$.

As a consequence, we obtain a strongly polynomial algorithm to
solve the following problem: given a convex grid $G$ and a
function $h_0:B(G)\to\Zset$, decide whether
$h_0$ is extendable to a concave cocirculation in $G$, and if so,
find an integer concave cocirculation $h$ with $h\rest{B(G)}=h_0$.
This is because the problem of finding a concave cocirculation having
prescribed values on the boundary can be written as a linear program
of size $O(n)$.

\smallskip
In view of the relationship of concave cocirculations and honeycombs,
one can give an analog of Theorem~\ref{tm:main} for honeycombs; we
omit it here.

\smallskip
In fact, one can slightly modify the method of proof of
Theorem~\ref{tm:main} so as to obtain the following strengthening: for
a concave cocirculation $h$ in a convex grid $G$, there exists an
integer concave cocirculation $h'$ satisfying $h'(e)=h(e)$ for each
edge $e\in O_h\cup I'_h$, where $I'_h$ is the set of edges contained
in circuits $C$ of $G$ such that $h$ takes integer values on all edges
of edges $C$. We omit the proof here.

\smallskip
Next, as mentioned in the Introduction, a cocirculation $h$ in
a convex grid $G$ need not admit an improvement to an integer
concave cocirculation in $G$ preserving the values on {\em all} edges
where $h$ is integral. (Note that in our proof of
Theorem~\ref{tm:main}, a vertex of the original honeycomb having only
one integer dual coordinate may split into vertices not preserving
this coordinate.) A counterexample $(G,h)$ is shown in the picture
(the right figure illustrates the corresponding honeycomb; here all
edge weights are ones, the integral edges are drawn in bold, and each
of the vertices $u,v,z$ has one integer dual coordinate).
 \begin{center}
  \unitlength=1mm
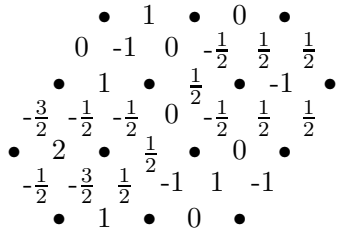
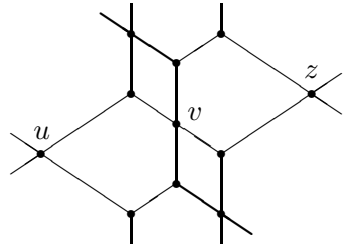
  \begin{picture}(120,32)
\put(6,1){\circle*{1.5}}
\put(18,1){\circle*{1.5}}
\put(30,1){\circle*{1.5}}
\put(0,10){\circle*{1.5}}
\put(12,10){\circle*{1.5}}
\put(24,10){\circle*{1.5}}
\put(36,10){\circle*{1.5}}
\put(6,19){\circle*{1.5}}
\put(18,19){\circle*{1.5}}
\put(30,19){\circle*{1.5}}
\put(42,19){\circle*{1.5}}
\put(12,28){\circle*{1.5}}
\put(24,28){\circle*{1.5}}
\put(36,28){\circle*{1.5}}
\put(11,0){1}
\put(23,0){0}
\put(1.2,4.5){-$\frac12$}
\put(7.2,4.5){-$\frac32$}
\put(13.5,4.5){$\frac12$}
\put(19.5,4.8){-1}
\put(26,4.8){1}
\put(31.5,4.8){-1}
\put(5,9){2}
\put(17,8.7){$\frac12$}
\put(29,9){0}
\put(1.2,13.5){-$\frac32$}
\put(7.2,13.5){-$\frac12$}
\put(13.2,13.5){-$\frac12$}
\put(20,13.8){0}
\put(25.2,13.5){-$\frac12$}
\put(32,13.5){$\frac12$}
\put(38,13.5){$\frac12$}
\put(11,18){1}
\put(23,17.7){$\frac12$}
\put(34,18){-1}
\put(8,22.8){0}
\put(13.2,22.8){-1}
\put(20,22.8){0}
\put(25.2,22.5){-$\frac12$}
\put(32,22.5){$\frac12$}
\put(38,22.5){$\frac12$}
\put(17,27){1}
\put(29,27){0}
\put(92,4){\circle*{1.0}}
\put(104,4){\circle*{1.0}}
\put(98,8){\circle*{1.0}}
\put(80,12){\circle*{1.0}}
\put(104,12){\circle*{1.0}}
\put(98,16){\circle*{1.0}}
\put(92,20){\circle*{1.0}}
\put(116,20){\circle*{1.0}}
\put(98,24){\circle*{1.0}}
\put(92,28){\circle*{1.0}}
\put(104,28){\circle*{1.0}}
\put(92,4){\line(-3,2){16}}
\put(92,4){\line(3,2){6}}
\put(92,20){\line(-3,-2){16}}
\put(92,20){\line(3,-2){12}}
\put(98,8){\line(0,1){16}}
\put(98,24){\line(3,2){6}}
\put(104,12){\line(3,2){16}}
\put(104,28){\line(3,-2){16}}
\put(79,14){$u$}
\put(99.5,16.5){$v$}
\put(115,22){$z$}
\thicklines
\put(92,4){\line(0,-1){4}}
\put(92,20){\line(0,1){12}}
\put(98,24){\line(-3,2){10}}
\put(98,8){\line(3,-2){10}}
\put(104,12){\line(0,-1){12}}
\put(104,28){\line(0,1){4}}
  \end{picture}
 \end{center}
One can check that $h$ is determined uniquely by its integer values,
i.e., $\Cscr(G,h\rest{F})=\{h\}$, where
$F$ is the set of edges where $h$ is integral.

\smallskip
We finish with the following question motivated by some aspects in
Section~\ref{sec:frac}. For a convex grid $G$, can one give a
``combinatorial characterization'' for the set of tilings of concave
cocirculations $h$ such that $h$ is a vertex of the polytope
$\Cscr(G,h\rest{B(G)})$?

\smallskip
{\bf Acknowledgement.} I thank Vladimir Danilov for stimulating
discussions and Gleb Koshevoy for pointing out to me a polyhedral
result in~\cite{DM}.

\end{document}